\documentclass[10pt]{amsart}

\usepackage{amsmath, amsfonts, amssymb, amsthm, mathdots, mathrsfs, dsfont}

\usepackage{tikz}

\usepackage[colorlinks, urlcolor=black]{hyperref}
\hypersetup{citecolor=blue, linkcolor=blue}

\newtheorem{theorem}{Theorem}[section]
\newtheorem{definition}[theorem]{Definition}
\newtheorem{proposition}[theorem]{Proposition}
\newtheorem{conjecture}[theorem]{Conjecture}

\newtheorem{lemma}[theorem]{Lemma}
\newtheorem{remark}[theorem]{Remark}

\numberwithin{equation}{section}

\author{Luis Alberto Lomel\'i}
\address{Luis Alberto Lomel\'i \\ Instituto de Matem\'aticas \\ Pontificia Universidad Cat\'olica de Valpara\'iso \\ Blanco Viel 596 \\ Cerro Bar\'on \\ Valpara\'iso \\ Chile}
\email{luis.lomeli@pucv.cl}

\begin{document}

\title[Langlands Program and Ramanujan]{Langlands Program and \\ Ramanujan Conjecture: a survey}

\begin{abstract}
We present topics in the Langlands Program to graduate students and a wider mathematically mature audience. We study both global and local aspects in characteristic zero as well as characteristic $p$. We look at modern approaches to the generalized Ramanujan conjecture, which is an open problem over number fields, and present known cases over function fields. In particular, we study automorphic $L$-functions via the Langlands-Shahidi method. Hopefully, our approach to the Langlands Program can help guide the interested reader into an exciting field of mathematical research.
\end{abstract}

\maketitle

\section*{Introduction}

Langlands' functoriality conjectures encompass a vast generalization of classical reciprocity laws that provide crossroads between Number Theory and Representation Theory. We let $\bf G$ be a connected reductive group and $\Gamma$ a certain subgroup. With the aid of Einsenstein series, Langlands made the important connection between the theory of Automorphic Forms and Representation Theory. Namely, an automorphic representation occurs as an irreducible constituent of the right regular action on the $\mathscr{L}^2$-space of complex valued functions
\begin{equation*}
   \mathscr{L}^2(\Gamma \backslash G) = 
   \mathscr{L}^2(\Gamma \backslash G)_{\rm disc} \oplus 
   \mathscr{L}^2(\Gamma \backslash G)_{\rm cont},
\end{equation*}
decomposed into discrete and continuous spectra. Initially, the group $\bf G$ was taken to be defined over the real numbers $\mathbb{R}$, yet we are able to study the spectral decomposition in the general setting of the ring of ad\`eles over a global field $k$ \cite{La1976,MoWa1994}. That is, $k$ is either a finite extension of the rational numbers $\mathbb{Q}$ (a number filed) or a finite separable extension of the field of functions for the projective line $\mathbb{P}^1$ defined over the finite field $\mathbb{F}_q$ (a function field).

Ubiquitous in modern number theory is the case of $G = {\rm SL}_2(\mathbb{R})$, with its principal congruence subgroup $\Gamma = {\rm SL}_2(\mathbb{Z})$. Modular forms find themselves inside the discrete spectrum $\mathscr{L}^2(\Gamma \backslash G)_{\rm disc}$. Of particular importance in the theory is Ramanujan's weight 12 modular form
\begin{equation*}
   \Delta(z) = q \prod_{n=1}^\infty (1-q^n)^{24}
   = \sum_{n=1}^\infty \tau(n) q^n,
\end{equation*}
where $q = {\rm exp}(2\pi iz)$ and $z$ is a variable in the hyperbolic upper half plane $\mathfrak{h}$. The coefficients $\tau(n)$ appearing in the Fourier expansion of $\Delta(z)$ give Ramanujan's $\tau$-function. It is a multiplicative function, i.e., $\tau(mn) = \tau(m)\tau(n)$ if $(m,n)=1$, and it satisfies the recursive relation on primes:
\begin{equation*}
   \tau(p^{n+1}) = \tau(p) \tau(p^n) - p^{11} \tau(p^{n-1}).
\end{equation*}
The original Ramanujan Conjecture states that
\begin{equation*}
   \left| \tau(p) \right| \leq 2p^{\frac{11}{2}},
\end{equation*}
and was proved by Deligne as a consequence of his work on the Weil conjectures \cite{De1974}. He proved it for modular forms, which is also known as the Ramanujan-Petersson conjecture. The Ramanujan conjecture generalizes to automorphic representations, which we study in this article.

Given a a global field $k$, its ring of ad\`eles $\mathbb{A}_k$ has the correct topology in order for us to do Fourier analysis. Namely, it is a restricted direct product taken over all of the completions over the primes of $k$, in addition to the arquimedean valuations, giving us that $\mathbb{A}_k$ is a locally compact abelian group. If we let ${\bf G} = {\rm GL}_1(\mathbb{A}_k)$ and $\Gamma = k^\times$, an automorphic representation in this setting is a Gr\"o\ss encharakter $\chi : k^\times \backslash \mathbb{A}_k^\times \rightarrow \mathbb{C}^\times$. We then find ourselves in the abelian case of Tate's thesis \cite{Ta1950}. The abelian case in general provides the base where Langlands functoriality conjectures build upon, namely Class Field Theory \cite{ArTa}.

We observe that by taking $k=\mathbb{Q}$ and the trivial Gr\"o\ss encharakter, we retrieve the Riemann $\zeta$-function from an automorphic $L$-function for ${\rm GL}_1$. i.e,
\begin{equation*}
   L^S(s,\mathds{1}) = \prod_p \left( 1 -p^{-s} \right)^{-1}
   											= \sum_{n=1}^\infty n^{-s} = \zeta(s),
\end{equation*}
where $S = \left\{ \infty \right\}$, and we initially have absolute convergence for $\Re(s) >1$. Tate $L$-functions were extended to automorphic $L$-functions for ${\rm GL}_n$ by Godement and Jacquet \cite{GoJa1972}, via an integral representation. Local Rankin-Selberg producs for representations of ${\rm GL}_m$ and ${\rm GL}_n$ were obtained by Jacquet, Piatetski-Shapiro and Shalika in \cite{JaPSSh1983}; for remarks and further references on the global theory we refer to \cite{CoPS2004} and \cite{HeLo2013a}. After preliminaries in Number Theory and Representation Theory, we present a uniform treatment of these $L$-functions and products for globally generic representations of the classical groups via the Langlands-Shahidi method.

Now, let $G = {\rm GL}_n(\mathbb{A}_k)$ and $\Gamma = {\bf Z}(\mathbb{A}_k) {\rm GL}_n(k)$, where $\bf Z$ denotes the center of ${\rm GL}_n$. The group ${\rm GL}_n$ is self dual, i.e., the connected component of its Langlands dual group is ${}^L{\rm GL}_n^\circ = {\rm GL}_n(\mathbb{C})$. The study of automorphic forms in $\mathscr{L}^2(\Gamma \backslash G)$ leads us to an important non-abelian generalized reciprocity law. The global Langlands conjecture states that there is a correspondence
\begin{equation*}
   \left\{ \begin{array}{c} \Pi \text{ a cuspidal automorphic} \\
  				      \text{representation of } {\rm GL}_n(\mathbb{A}_k)
					\end{array} \right\}
   \longleftrightarrow 
   \left\{ \begin{array}{c} \Phi_\Pi: \mathcal{L}_k \rightarrow {}^L{\rm GL}_n  \text{ an } \\
   				     \text{irreducible $L$-parameter}
					\end{array} \right\}
\end{equation*}
If $k$ is a global function field, the group $\mathcal{L}_k$ is essentially the Weil group $\mathcal{W}_k$ and we note that global Langlands is a landmark result of V. G. Drinfeld for ${\rm GL}_2$ \cite{Dr1978,Dr1988} and of L. Lafforgue for ${\rm GL}_n$ \cite{LaL2002}. The case of a number field $k$ being open, there should be a conjectural locally compact topological group $\mathcal{L}_k$ associated to the absolute Galois group $\mathcal{G}_k$, where $\mathcal{L}_k$ should have the right topology and there should be a morphism $\mathcal{L}_k \rightarrow \mathcal{G}_k$ for us to have a one-to-one correspondence between automorphic representations of ${\rm GL}_n$ and Galois representations, as mentioned by Langlands in \cite{LaCorvallisII}. The correspondence is in such a way that it respects global $L$-functions and root numbers, and agrees with the local Langlands correspondence at every place. On the left hand side we have automorphic $L$-functions and on the right hand side we have Artin $L$-functions.
The correspondence is known to imply the Ramanujan Conjecture for automorphic representations.

%(The Weil group $\mathcal{W}_k$ is a locally pro-finite topological group closely related to the absolute Galois group $\mathcal{G}_k$ of $k$, which has the right topology for the study of known cases of Langlands correspondences. In the function field case, we study $\ell$-adic representations of $\mathcal{W}_k$ and over a non-archimedean local field $F$ we study Weil-Deligne representations of $\mathcal{W}_F$.)

When we take the completion of a global field $k$ at one of its absolute values $v$, we obtain a local field, and all local fields arise in this way. A non-archimedean local field is obtained from a $\mathfrak{p}$-adic valuation, $\mathfrak{p}$ a prime ideal in the ring of integers $\mathcal{O}_k$ of a global field $k$. We let $F$ denote a non-archimedean local field. 

We work with connected quasi-split reductive groups over local and global fields. If the reader is unfamiliar with reductive groups, it is safe to assume on a first reading of this article that $\bf G$ is a group given by the examples of section \S~\ref{examples}, which include the classical groups. Given a connected reductive group $\bf G$ defined over the local field $F$, we write $G = {\bf G}(F)$. We study admissible representations $\pi$ of $G$ over a complex infinite dimensional vector space $V$; we refer to the first three sections of \cite{BuHe2006} for basic notions of admissible representations and further references.

The local Langlands correspondence for ${\rm GL}_n$ provides us with a one to one reciprocity relation
\begin{equation*}
   \left\{ \begin{array}{c} \pi \text{ an irreducible admissible} \\
  				      \text{representation of } {\rm GL}_n(F)
					\end{array} \right\}
   \longleftrightarrow 
   \left\{ \begin{array}{c} \phi_\pi: \mathcal{W}_F' \rightarrow {}^L{\rm GL}_n \\
   				     \text{an } L\text{-parameter}
					\end{array} \right\},
\end{equation*}
%\begin{equation*}
   %\left\{ \begin{array}{c} \pi \text{ an irreducible admissible} \\
  	%			      \text{representation of } {\rm GL}_n(F)
	%				\end{array} \right\}
   %\longleftrightarrow 
   %\left\{ \begin{array}{c} \sigma \text{ an $n$-dimensional ${\rm Frob}$-semisimple} \\
   	%			     \text{Weil-Deligne representation of } \mathcal{W}_F
	%				\end{array} \right\},
%\end{equation*}
where the Weil group $\mathcal{W}_F$ provides us with the right locally compact group on the Galois side, yet we need $L$-parameters of the Weil-Deligne group $\mathcal{W}_F'$ or Weil-Deligne representations of $\mathcal{W}_F$, see \S~7 of \cite{BuHe2006}. The correspondence is uniquely characterized by the fact that it preserves local $L$-functions and root numbers for products of representations \cite{He1993}.

From the theory of admissible representations of a $\mathfrak{p}$-adic group $G = {\bf G}(F)$, an irreducible admissible representation $\pi$ is supercuspidal if it cannot be obtained as a constituent of a parabolically induced representation
\begin{equation*}
   {\rm Ind}_P^G (\rho).
\end{equation*}
Here, $\bf P$ is a parabolic subgroup of $\bf G$, with Levi component $\bf M$ and unipotent radical $\bf N$; $\rho$ is a representation of $M$, extended trivially on $N$; and, we use unitary parabolic induction. For ${\bf G} = {\rm GL}_n$, the Levi subgroups are all of the form $\prod {\rm GL}_{n_i}$, with $n =\sum n_i$, and standard parabolic subgroups are block upper triangular. The Levi subgroups and parabolic subgroups of the classical groups are described in \S~\ref{examples}. As observed in \cite{BuKu1993}, there are two main steps in undesrtading admissible representations: (i) constructing supercuspidal representations; (ii) constructing admissible representations via parabolic induction.

The first examples of supercuspidals were obtained in \cite{Ma1964}, and work using the theory of types of \cite{BuKu1993} gives all supercuspidals for ${\rm GL}_n$. Stevens constructs all supercuspidals for the classical groups, assuming ${\rm char}(F) \neq 2$ \cite{St2008}. For a general reductive group, the approach of Yu gives tame supercuspidals \cite{Yu2001,KiJL2007}. And we have the supercuspidals of Gross and Reeder, in addition to the work of Kaletta \cite{Ka2016}. The Kazhdan transfer between close local fields provides an approach to obtain the local Langlands correspondence in characteristic $p$ from known cases in characteristic $0$ \cite{AuBaPlSo2016, GaVa2017}.

We observe that studying the space $\mathscr{L}^2(G)$ leads to the Plancherel Theorem much studied by Harish-Chandra \cite{Wa2003}. For a classical group $\bf G$, in addition to $L$-functions and root numbers, we use Harish-Chandra $\mu$-functions in order to provide a characterization of the local Langlands transfer. An important class of representations of a $\mathfrak{p}$-adic reductive groups are tempered representations; we recall that an admissible representation is tempered if its matrix coefficients are in $\mathscr{L}^{2+\varepsilon}(G)$ for every $\varepsilon >0$.

\smallskip

\noindent{\bf Ramanujan Conjecture for ${\rm GL}_n$.} \emph{Let $\Pi$ be a cuspidal automorphic representation of ${\rm GL}_n(\mathbb{A}_k)$, then each local component $\Pi_v$ is tempered.
}

\smallskip

If $\Pi_v$ is unramified, temperedness implies that its Satake parameters satisfy
\begin{equation*}
   \left| \alpha_{j,v} \right|_{k_v} = 1, \quad 1 \leq j \leq n.
\end{equation*}
A priori, an arbitrary constituent $\Pi_v$ of a cuspidal automorphic representation is admissible, although it is unramified for almost every $v$. Langlands' classification allows us to express general admissible representations in terms of tempered data \cite{BoWa, Si1978}. In order to understand the bounds towards the Ramanujan conjecture, let us recall Langlands classification for ${\bf G}={\rm GL}_n$. Let $\nu$ be the character of a general linear group over $F$ given by $\nu(g) = \left| \det(g) \right|_F$. Then, if $\pi$ is an admissible representation of $G$, it has a unique Langlands quotient $\xi$ obtained from the induced representation
\begin{equation*}
   {\rm ind}_P^G(\pi_1 \otimes \cdots \otimes \pi_b).
\end{equation*}
Here, each $\pi_i$ is a representation of ${\rm GL}_{n_i}(F)$ of the form
\begin{equation*}
   \pi_j = \tau_j \cdot \nu^{r_j},
\end{equation*}
with each $\tau_j$ tempered, and we can assume that the Langlands parameters are arranged so that $0 \leq r_1 \leq \cdots \leq r_b$.

The Ramanujan conjecture for ${\rm GL}_n$ states that $r_b = 0$ for every component $\Pi_v$ of a cuspidal automorphic $\Pi$. If $k$ is a function field, the conjecture was established for ${\rm GL}_2$ by Drinfeld \cite{Dr1978} and ${\rm GL}_n$ by L. Lafforgue \cite{LaL2002}. For $\bf G$ different than ${\rm GL}_n$, we have a discussion in \S~\ref{GenRamanujan}, since it is no longer true for all cuspidal automorphic representations and it is rephrased for globally generic representations. The Ramanujan conjecture is a result of the author for the classical groups over function fields \cite{Lo2009,LoRationality}.

Let $k$ be a number field. In this open case, we have bounds towards Ramanujan. A first bound, known as the trivial bound is $r_b \leq 1/2$ \cite{JaSh1981}. The best known bounds for arbitrary $n$, were obtained by Luo, Rudnick and Sarnak in \cite{LuRuSa1999}, expressed by the following formula for the Ramanujan bound
\begin{equation*}
   r_b \leq \dfrac{1}{2} - \dfrac{1}{n^2+1}.
\end{equation*}
However, for small $n$, there are better results. First, for $n=2$, the best known bound of $r_b \leq 7/64$ was obtained for $k=\mathbb{Q}$ by Kim and Sarnak in \cite{KiSa2003}. For an arbitrary number field $k$, this Ramanujan bound was proved by Blomer and Brumley in \cite{BlBr2011}, where one can also find the best known bounds for up to $n=4$.

We mention the work by Clozel and Thorne, starting with \cite{ClTh2014}, which can lead towards newer Ramanujan bounds for ${\rm GL}_2$ over number fields. For ${\rm GL}_n$, certain cuspidal automorphic representations are known to satisfy Ramanujan over number fields, where there is an ongoing 10 author work towards potential automorphy \cite{ACC+}. Over function fields, there is the important work by V. Lafforgue over function fields, who establishes an automorphic to Galois functorial lift \cite{LaV}.

This survey article presents an approach to functoriality in a uniform way for a global field, either characteristic zero or $p$, via Langlands-Shahidi automorphic $L$-functions and the Converse Theorem of Cogdell and Piatetski-Shapiro. We do not touch upon the trace formula, where Arthur has established the Endoscopic classification for automorphic representations of the classical groups over number fields, relying on the stabilization of the twisted trace formula completed by M\oe glin and Waldspurger in \cite{MoWaX}. The Fundamental Lemma, proved by Ngo providing a most intricate crucial ingredient.

There are several survey articles available, many of them focusing on number fields. General introductions to the Langlands Program are given by Arthur \cite{Ar2003}, Gelbart \cite{Ge1984}, and a book that may be found useful by graduate students \cite{BeGe2004}. Frenkel touches upon the geometric Langlands Program \cite{Fr2013}. And, there is much recent progress in the mod-$\ell$ and the $p$-adic Langlands Program; a nice read and further references may be found in \cite{Au2018}. Furthermore, the reader may find excelent surveys for the Ramanujan Conjecture over number fields by Blomer-Brumley \cite{BlBr2013}, Sarnak \cite{Sa2005}, and Shahidi \cite{ShF2004}, for example.

\subsection*{Acknowledgments} The author thanks the organizers of the reunions: Mexican Mathematicians in the World 2016 and 2018; and the institutions CIMAT, Guanajuato, and Banff, Casa M\'exico Oaxaca, where the conferences were held; a rich atmosphere for academic exchange was provided and the idea for this article came to be. He is grateful to IISER, Pune, and TIFR, Mumbay, for an invitation to conduct research in January-February 2017, where a first draft of this article was written. He would like to thank A.-M. Aubert, G. Henniart, L. Lafforgue, V. Lafforgue, M. Mishra, D. Prasad, F. Shahidi, S. Varma for mathematical communications and helpful discussions. The author was partially supported by Project VRIEA-PUCV 039.367 and FONDECYT Grant 1171583.

\section{Artin $L$-functions and Number Theory}

Let $F$ be a non-archimedean local field. Thus, $F$ is either an extension of the $p$-adic numbers $\mathbb{Q}_p$ or it is of characteristic $p$ and of the form $F \simeq \mathbb{F}_q (\!\left( t \right)\!)$, the field of Laurent power series with coefficients in the finite field $\mathbb{F}_q$. We let $\mathcal{O}_F$ be the compact open ring of integers of $F$, it is a discrete valuation ring with a unique maximal compact open prime ideal $\mathfrak{p}_F = (\varpi_F)$, where $\varpi_F$ is a uniformizer. We refer to \S~4 of \cite{RaVa1999} for the classification of locally compact fields and their arithmetic properties.

\subsection{}\label{galois} The absolute Galois group
\begin{equation*}
   \mathcal{G}_F = \lim_{\longleftarrow} {\rm Gal}(E/F),
\end{equation*}
where $E/F$ varies through finite Galois extensions $E$ over $F$ inside a separable algebraic closure $\overline{F}$, is a pro-finite group with the Krull topology, hence compact. We obtain isomorphic Galois groups if we change the separable closure. Initially we would like to study representations of $\mathcal{G}_F$, however, we study the related locally pro-finite Weil group $\mathcal{W}_F$ in order to obtain the Langlands correspondence. In fact, we study isomorphism classes of representations, since this makes them independent of the choice of separable closure $\overline{F}$, which we fix.

Algebraically, the Weil group $\mathcal{W}_F$ is generated by the Frobenius elements and is normal in $\mathcal{G}_F$. The inertia group $\mathcal{I}_F = {\rm Gal}(\overline{F}/F_{\rm nr})$, where $F_{\rm nr}$ is the maximal unramified extension of $F$, is embedded inside $\mathcal{W}_F$. We denote the maximal tamely ramified extension of $F$ by $E_{\rm t}$, it is obtained as the composite of all the cyclic extensions $E_n$ of $F_{\rm nr}$. The wild inertia group is then $\mathcal{P}_F = {\rm Gal}(\overline{F}/E_{\rm t})$. The topology of $\mathcal{W}_F$ is such that the inertia group $\mathcal{I}_F$ is seen as an open subgroup of $\mathcal{W}_F$ and coincides with the topology of ${\rm Gal}(\overline{F}/F_{\rm nr})$ inside $\mathcal{G}_F$. We express local class field theory in terms of the Artin map in Theorem~\ref{local:artin}, and we now recall the process of ``abelianization'' and the ``Verlagerung'' homomorphism.

Let ${\rm Der}(\mathcal{W}_F)$ be the derived group and $\overline{{\rm Der}(\mathcal{W}_F)}$ its closure, then $\mathcal{W}_F^{\rm ab} = \mathcal{W}_F / \overline{{\rm Der}(\mathcal{W}_F)}$ is the abelianization of $\mathcal{W}_F$. Given a separable extension of non-archimedean local fields $E/F$, we have the norm map ${\rm N}_{E/F}: E^\times \rightarrow F^\times$ and
\begin{equation*}%\label{Verlagerung}
   {\rm Ver}_{E/F}: \mathcal{W}_F^{\rm ab} \rightarrow \mathcal{W}_E^{\rm ab}
\end{equation*}
is the homomorphism defined by its action on cosets
\begin{equation*}
   x \, {\rm Der}(\mathcal{W}_F) \quad \mapsto \prod_{y \, \in \mathcal{W}_F / \mathcal{W}_E} y \, {\rm Der}(\mathcal{W}_E)
\end{equation*}
and passing to the closure on derived groups $\overline{{\rm Der}(\mathcal{W}_F)} \subset \mathcal{W}_F^{\rm ab}$ and $\overline{{\rm Der}(\mathcal{W}_E)} \subset \mathcal{W}_E^{\rm ab}$.

Weil groups were introduced in \cite{We1951}. An excellent reference for local Weil groups and their representations can be found in chapter~7 of \cite{BuHe2006}. We will turn towards the global theory, which is existential in nature, in \S~\ref{Weil:global}.

\subsection{Local class field theory}\label{local:artin}
\emph{For every non-archimedean field $F$ there is a reciprocity map
\begin{equation*}
   \mathcal{W}_F \rightarrow {\rm GL}_1(F),
\end{equation*}
sending $\mathcal{I}_F$ to $\mathcal{O}_F^\times = {\rm GL}_1(\mathcal{O}_F)$, $\mathcal{P}_F$ to $1+\mathfrak{p}_F$, and giving a topological isomorphism 
\begin{equation*}
   {\rm Art} : {\rm GL}_1(F) \rightarrow \mathcal{W}_F^{\rm ab},
\end{equation*}
between the abelianization $\mathcal{W}_F^{\rm ab}$ of $\mathcal{W}_F$ and $F^\times = {\rm GL}_1(F)$. The Artin map is compatible with finite separable extensions, for which we have commuting diagrams}
\begin{center}
\begin{tikzpicture}
%\tikzstyle{suite}=[rounded corners]
   \draw (0,0) node {$E^\times$};
   \draw[<-] (0,.25) -- (0,1);
   %\draw (0,.625) node [left] {$\alpha$};
   \draw (0,1.25) node {$\mathcal{W}_E$};
   \draw[->] (.5,1.25) -- (1.5,1.25);
   \draw (.5,1.25) arc (270:90:.0625);
   %\draw (1,1.25) node [above] {$\subset$};
   \draw (2,1.25) node {$\mathcal{W}_F$};
   \draw (2,0) node {$F^\times$};
   \draw[->] (.5,0) -- (1.5,0);
   \draw (1,0) node [below] {${\rm N}_{E/F}$};
   \draw[->] (2,1) -- (2,.25);
   \draw (3,.625) node {and};
   \draw (6,1.25) node {$\mathcal{W}_E^{\rm ab}$};
   \draw (6,0) node {$E^\times$};
   \draw[->] (4.5,0) -- (5.5,0);
   \draw (4.5,0) arc (270:90:.0625);
   \draw[->] (6,1) -- (6,.25);
   %\draw (6,.625) node [right] {$\alpha$};
   \draw (4,0) node {$F^\times$};
   \draw[<-] (4,.25) -- (4,1);
   %\draw (4,.625) node [left] {$\alpha$};
   \draw (4,1.25) node {$\mathcal{W}_F^{\rm ab}$};
   \draw[->] (4.5,1.25) -- (5.5,1.25);
   \draw (5,1.25) node [above] {${\rm Ver}_{E/F}$};
\end{tikzpicture}
\end{center}

\subsection{Global fields and valuations} Let $k$ be a global field. Thus, $k$ is either a finite algebraic extension of $\mathbb{Q}$, i.e. a number field, or a finite separable extension of $\mathbb{F}_q(t)$, i.e. a function field. Let $\mathcal{O}_k$ be the ring of integers of $k$, it is a Dedekind domain.

A valuation $v$ of $k$ can be either non-archimedean or archimedean, the latter possible only in the number field case. We write ${\rm Val}(k)$ for the set of valuations of $k$. After localization, we obtain a local field $k_v$, where $k_v \simeq \mathbb{R}$ or $\mathbb{C}$ if $v \in {\rm Val}(k)$ is archimedean. We write ${\rm Val}_\mathfrak{p}(k)$ for the set of non-archimedean valuations of $k$. Every $v \in {\rm Val}_\mathfrak{p}(k)$ thus leads to a non-archimedean local field $k_v$, with corresponding ring of integers $\mathcal{O}_v$ and maximal prime ideal $\mathfrak{p}_v$. We write $\mathfrak{p}_v = \left( \varpi_v \right)$ for the maximal prime ideal, $\varpi_v$ a uniformizer, and $q_v$ the cardinality of the finite residue field $\mathbb{F}_v = \mathcal{O}_v / \mathfrak{p}_v$.

Let $l/k$ be a finite separable extension, so we have an extension of Dedekind rings $\mathcal{O}_k \subset \mathcal{O}_l$. Given $v \in {\rm Val}_\mathfrak{p}(k)$, we have that $\mathfrak{p}_v$ factorizes in $\mathcal{O}_l$ in such a way that one obtains all valuations above $v$
\begin{center}
\begin{tikzpicture}
   \draw (0,.625) node {$\mathfrak{p}_v \mathcal{O}_l = \mathfrak{P}_{w_1}^{e_1} \cdots \mathfrak{P}_{w_g}^{e_g}$};
   \draw[<->] (2.25,.625) -- (3.25,.625);
   \draw (4,1.25) node {$w_1$};
   \draw (5,1.25) node {$\cdots$};
   \draw (6,1.25) node {$w_g$};
   \draw[-] (5,0) -- (4,1);
   \draw (5,0) node [below] {$v$};
   \draw[-] (5,0) -- (6,1);
   \draw (7.5,0) node {${\rm Val}_\mathfrak{p}(k)$};
   \draw (7.5,1.25) node {${\rm Val}_\mathfrak{p}(l)$};
\end{tikzpicture}
\end{center}
where $w_i \longleftrightarrow \mathfrak{P}_i$, for $1 \leq i \leq g$. We write $w \vert v$ when $w=w_i$ for some $i$ as above.

Each positive integer $e_w$ is called the ramification index of $l_{w} / k_v$, and we also have $f_w$ such that $q_{w} = q_v^{f_w}$. We say that $v \in {\rm Val}_\mathfrak{p}(k)$ is unramified if $e_w = 1$ for every $w \vert v$. It is a known fact that almost every $v \in {\rm Val}_\mathfrak{p}(k)$ is unramified, i.e., all but finitely many. We further have
\begin{equation*}
   n = [l:k] = \sum_{w \vert v} [l_w:k_v] = \sum_{w \vert v} e_w f_w.
\end{equation*}
If $l/k$ is Galois, then $e_w = e$ and $f_w = f$ for all $w \vert v$ and we have $n = efg$. For more details see, for example, \S~4.4 of \cite{RaVa1999}.

\subsection{On Artin $L$-functions} Fix a Galois extension $l/k$. The Galois group ${\rm Gal}(l/k)$ transitively permutes the places $w \vert v$ and each $\alpha \in {\rm Gal}(l/k)$ induces an isomorphism
\begin{equation*}
   \alpha_w : l_w \simeq l_{\alpha(w)}.
\end{equation*}
The decomposition group is defined by
\begin{equation*}
   G_w = \left\{ \alpha \in {\rm Gal}(l/k) \ \vert \ \alpha(w) = w \right\},
\end{equation*}
for each $w \vert v$. We have that
\begin{equation*}
   G_{\beta(w)} = \beta \, G_w \beta^{-1}.
\end{equation*}
If $\alpha \in G_w$, then $\alpha$ induces a $k_v$-automorphism $\alpha_w$ on $l_w$, and we have
\begin{equation*}
   G_w \xrightarrow{\sim} {\rm Gal}(l_w/k_v) \twoheadrightarrow 
   {\rm Gal}(\mathbb{F}_w/\mathbb{F}_v).
\end{equation*}
The kernel of the last surjective homomorphism is the inertia subgroup $\mathcal{I}_w \subset G_w \subset {\rm Gal}(l/k)$. And for almost all $v \in {\rm Val}(k)$, we have that $\mathcal{I}_w$ is trivial. Indeed, at unramified places we have
\begin{equation*}
   {\rm Gal}(l_w/k_v) \simeq {\rm Gal}(\mathbb{F}_w/\mathbb{F}_v),
\end{equation*}
which is cyclic and generated by
\begin{equation*}
   \bar{\alpha}_w: x \mapsto x^{q_v}.
\end{equation*}

Now, a Frobenius automorphism for $l/k$ is any element $\alpha_w \in G_w$ with image $\bar{\alpha}_w$ in ${\rm Gal}(\mathbb{F}_w/\mathbb{F}_v)$. It is characterized by the property
\begin{equation*}
   \alpha_w(x) = x^{q_v} \, ({\rm mod} \, \mathfrak{P}_w), \ x \in \mathcal{O}_w.
\end{equation*}
We observe that $\alpha_{\beta(w)} = \beta^{-1} \alpha_w \beta$, and hence $\alpha_w$ is determined up to conjugacy. The conjugacy class is by definition the Frobenius element of $l/k$ at the place $v$ ${\rm Frob}_{l/k,v}$.

We look at the global Weil group $\mathcal{W}_k$ in \S~\ref{Weil:global}, which for every Galois extension $l/k$ has an isomorphism
\begin{equation*}
   \mathcal{W}_k/\mathcal{W}_l \simeq {\rm Gal}(l/k).
\end{equation*}
Assume that we have an $n$-dimensional continuous representation $(\sigma,V)$ of $\mathcal{W}_k$ which factors through $\mathcal{W}_l$, for $l/k$ Galois. And consider $\sigma$ as a representation of ${\rm Gal}(l/k)$. If we let
\begin{equation*}
   V_w = V^{\mathcal{I}_w},
\end{equation*}
then $\sigma({\rm Frob}_v)$ is defined on $V_w$. In this setting, ``incomplete'' Artin $L$-functions are given by
\begin{equation*}
   L(s,\sigma) = \prod_{v \in {\rm Val}_\mathfrak{p}(k)} \dfrac{1}{\det(I - \sigma({\rm Frob}_v)\vert_{V_w} q_v^{-s})},
\end{equation*}
originally absolutely and uniformly convergent for $\Re(s) \geq 1+\delta$, $\delta > 0$.

\subsection*{Artin Conjecture} \emph{Artin $L$-functions $L(s,\sigma)$ are entire for non-trivial Galois representations of dimension $n>1$.}

\smallskip

The conjecture is expected in the setting of \S~\ref{Weil:global}, where Deligne and Langlands \cite{De1973,La1976} define local root numbers $\varepsilon(s,\sigma_v,\psi_v)$ at every place $v$ of $k$. Thus, we have a global Artin root number
\begin{equation*}
   \varepsilon(s,\sigma) = \prod_{v} \varepsilon(s,\sigma_v,\psi_v),
\end{equation*}
where $\varepsilon(s,\sigma_v,\psi_v) = 1$ at every unramified place. Artin $L$-functions can be extended to archimedean places $v$ of $k$, in a way that they are inductive in nature, and the resulting ``completed'' Artin $L$-function
\begin{equation*}
  L(s,\sigma) = \prod_{v \in {\rm Val}(k)} L_v(s,\sigma),
\end{equation*}
is in agreement with the previous definition
\begin{equation*}
   L_v(s,\sigma) = \dfrac{1}{\det(I - \sigma({\rm Frob}_v)\vert_{V_w} q_v^{-s})},
\end{equation*}
for $v \in {\rm Val}_\mathfrak{p}(k)$ and $\sigma$ factors through $\mathcal{W}_l$ as above. Artin $L$-functions in general have a meromorphic continuation to the complex plane. With the aid of Brauer theory, they can be proved to satisfy a functional equation
\begin{equation*}
   L(s,\sigma) = \varepsilon(s,\sigma) L(1-s,\tilde{\sigma}),
\end{equation*}
where $\tilde{\sigma}$ denotes the contragredient representation. This is done by reducing to the case of $n = 1$, studied by Artin and Tate. We discuss later on in \S~\ref{globalization} cases where the Artin Conjecture is known.

\subsection{Ad\`eles and Id\`eles}
Tate's celebrated thesis \cite{Ta1950} uses restricted direct products to make a study of Fourier Analysis using valuation vectors, as Chevalley called them, or in modern language the ad\`eles and id\`eles. 

Restricted direct products are defined over a global field $k$. We assume that we are given a a finite set $S_0$ consisting of valuations $v \in {\rm Val}(k)$, including all archimedean $v$, and a locally compact topological group $G_v$ for every $v \in {\rm Val}(k)$ such that $G_v$ has a compact open subgroup $H_v$ for every $v \notin S_0$. In this situation, we say that a condition is valid for almost all $v$ if there exists a finite set $S \supset S_0$ such that the condition holds for every $v \notin S$. We define the restricted direct product
\begin{equation*}
   G = \prod{}' (G_v:H_v)
\end{equation*}
as the set of elements of the form $(g_v)$, $g_v \in G_v$, where we have that $g_v \in H_v$ for almost all $v$.

The open sets of the group $G$ are given by the following basis of neighborhoods of the identity
\begin{equation*}
   U_S = \prod_{v \in S} U_v \times \prod_{v \notin S} H_v,
\end{equation*}
where $S \supset S_0$ is a finite set of places and, for each $v \in S$ we have an arbitrary open set $U_v$ in $G_v$. We observe that the resulting topology on $G$ makes it a locally compact topological group.

The ad\`eles $\mathbb{A}_k$ of a global field $k$ are defined by taking $S_0 = {\rm Val}(k) \setminus {\rm Val}_\mathfrak{p}(k)$ (the set of archimedean valuations), setting $G_v = k_v$ for every $v \in {\rm Val}(k)$ and $H_v = \mathcal{O}_v$ for $v \notin S_0$. The id\`eles are the group ${\rm GL}_1(\mathbb{A}_k)$, and it is not a problem for us define ${\rm GL}_n(\mathbb{A}_k)$, for any positive integer $n$. Indeed, observe that for every $v \in {\rm Val}_\mathfrak{p}(k)$, the group ${\rm GL}_n(\mathcal{O}_v)$ is a maximal compact open subgroup of ${\rm GL}_n(k_v)$. Then, ${\rm GL}_n(\mathbb{A}_k)$ is given by the restricted direct product
\begin{equation*}
   {\rm GL}_n(\mathbb{A}_k) = \prod {}' \left( {\rm GL}_n(k_v) : {\rm GL}_n(\mathcal{O}_v) \right).
\end{equation*}
Algebraically $\mathbb{A}_k^\times$ is the multiplicative group of $\mathbb{A}_k$, yet the topology of ${\rm GL}_1(\mathbb{A}_k)$ is not the one inherited from $\mathbb{A}_k$ by restriction. Nonetheless, we abuse notation and denote the id\`eles by $\mathbb{A}_k^\times$, and we make the convention that the topology is the one obtained as a restricted direct product $\mathbb{A}_k^\times = {\rm GL}_1(\mathbb{A}_k)$.

Tate's thesis can be seen as an analytic approach to the abelian case studied by Class Field Theory. Let $\psi = \otimes \psi_v$ be a non-trivial additive character of $k \backslash \mathbb{A}_k$, allowing us to do Fourier Analysis, and let
\begin{equation*}
   \chi: k^\times \backslash \mathbb{A}_k^\times \rightarrow \mathbb{C}^\times
\end{equation*}
be a Gr\"o\ss encharakter of $\mathbb{A}_k^\times = {\rm GL}_1(\mathbb{A}_k)$. In \cite{Ta1950}, Tate works over number fields and defines $L$-functions and root numbers
\begin{equation*}
   L(s,\chi) = \prod_v L(s,\chi_v) \text{ and } 
   \varepsilon(s,\chi) = \prod_v \varepsilon(s,\chi_v,\psi_v),
\end{equation*}
and establishes their functional equation
\begin{equation*}
   L(s,\chi) = \varepsilon(s,\chi) L(1-s,\chi^{-1}).
\end{equation*}
Tate's thesis can be studied over function fields as well, and a great reference valid over any global field $k$ can be found in \cite{RaVa1999}.

Godement and Jacquet extend Tate's theory of $L$-functions for Gr\"o\ss encharakters $\chi$ of ${\rm GL}_1(\mathbb{A}_k)$ to automorphic representations $\Pi$ of ${\rm GL}_n(\mathbb{A}_k)$ in \cite{GoJa1972}. We thus have global Godement-Jacquet $L$-functions and root numbers
\begin{equation*}
   L(s,\Pi) = \prod_v L(s,\Pi_v) \text{ and } 
   \varepsilon(s,\Pi) = \prod_v \varepsilon(s,\Pi_v,\psi_v),
\end{equation*}
which satisfy the functional equation
\begin{equation*}
   L(s,\Pi) = \varepsilon(s,\Pi) L(1-s,\tilde{\Pi}).
\end{equation*}
We turn towards automorphic representations in \S~\ref{rgrt} and we present all the automorphic $L$-functions that we will need in this article in a uniform way in \S~\ref{LSmethod}.

\subsection{Weil groups and Artin $L$-functions}\label{Weil:global}
Let $k$ be a global field and let $\mathcal{G}_k$ be its absolute Galois group obtained just as in \eqref{galois} after fixing a separable algebraic closure $\overline{k}$. The global Weil group is a topological group $\mathcal{W}_k$ together with a dense continuous homomorphism $\mathcal{W}_k \rightarrow \mathcal{G}_k$, such that for every finite separable extension $l/k$: $\mathcal{W}_l$ is open in $\mathcal{W}_k$ and maps to $\mathcal{G}_l$ with dense image; there is an isomorphism \begin{equation*}
   \mathcal{W}_l^{\rm ab} \simeq \mathcal{G}_l^{\rm ab}
\end{equation*} 
which is compatible with global abelian class field theory, see \S~\ref{global:artin}; and, the projective limit
\begin{equation*}
   \mathcal{W}_k = \lim_{\longleftarrow} (\mathcal{W}_k / \, \overline{\mathcal{W}_l})
\end{equation*}
taken over all separable $l$ inside $\overline{k}$ is an isomorphism of topological groups. We observe that the constructions of \S~1.1, e.g. abelianization and Verlagerung map, are valid in this setting. Furthermore, these constructions also respect isomorphisms of separable field extensions.

In order to take a closer look at Artin $L$-functions, write $\mathcal{W}_v = \mathcal{W}_{k_v}$, for every $v \in {\rm Val}_\mathfrak{p}(k)$. Since we mod out by the inertia group $\mathcal{I}_{v}$, we can define ${\rm Frob}_v$ on $\mathcal{W}_v$ from the isomorphism
\begin{equation*}
   {\rm Gal}(k_{v,\infty}/k_v) \cong \lim_{\longleftarrow} \, \mathbb{Z}/n\mathbb{Z} = \widehat{\mathbb{Z}}.
\end{equation*}
For each $n$, there is a unique unramified extension $k_{v,n}$ of $k_v$ with cyclic Galois group ${\rm Gal}(k_{v,n}/k_v)$ and $k_{v,\infty}$ is the composite of all $k_{v,n}$.

Artin $L$-functions at archimedean places $v$ of $k$ can be defined \cite{BeGe2004}, and this leads us to ``completed'' Artin $L$-functions. In \cite{De1973}, Deligne provides a proof that Artin $L$-functions satisfy the functional equation
\begin{equation*}
   L(s,\sigma) = \varepsilon(s,\sigma) L(1-s,\tilde{\sigma}),
\end{equation*}
where $\tilde{\sigma}$ denotes the contragredient representation.

The class field theory of Artin-Takagi, gives a very complete description of abelian Galois extensions in connection to the id\`eles. We next write it in the context of Weil groups, which are constructed cohomologically in Artin-Tate \cite{ArTa}, and we also note the important references \cite{CaFr} and \cite{Ta1979}.

\subsection{Global class field theory}\label{global:artin}
\emph{Let $k$ be a global field. For every finite separable extension $l/k$ there is a reciprocity map $\mathcal{W}_l \rightarrow {\rm GL}_1(\mathbb{A}_k)$ inducing an isomorphism
\begin{equation*}
   {\rm Art} : \mathcal{C}_l = l^\times \backslash {\rm GL}_1(\mathbb{A}_l) \xrightarrow{\sim} \mathcal{W}_l^{\rm ab}.
\end{equation*}
For abelian Galois extensions $l/k$ it induces an isomorphism
\begin{equation*}
   {\rm Art} : \Gamma \backslash {\rm GL}_1(\mathbb{A}_l)
   \xrightarrow{\sim} {\rm Gal}(l/k),
\end{equation*}
for $\Gamma = l^\times {\rm N}_{l/k} ({\rm GL}_1(\mathbb{A}_l))$, where ${\rm N}_{l/k}$ is the global norm. The Artin map is compatible with finite separable extensions, for which we have commuting diagrams
\begin{center}
\begin{tikzpicture}
%\tikzstyle{suite}=[rounded corners]
   \draw (0,0) node {$\mathcal{W}_l^{\rm ab}$};
   \draw[<-] (0,.25) -- (0,1);
   %\draw (0,.625) node [left] {$\alpha$};
   \draw (0,1.25) node {$\mathcal{C}_l$};
   \draw (1,1.25) node [above] {${\rm N}_{l/k}$};
   \draw[->] (.5,1.25) -- (1.5,1.25);
   \draw (.5,0) arc (270:90:.0625);
   %\draw (1,1.25) node [above] {$\subset$};
   \draw (2,1.25) node {$\mathcal{C}_k$};
   \draw (2,0) node {$\mathcal{W}_k^{\rm ab}$};
   \draw[->] (.5,0) -- (1.5,0);
   \draw[->] (2,1) -- (2,.25);
   \draw (3,.625) node {and};
   \draw (6,1.25) node {$\mathcal{C}_l$};
   \draw (6,0) node {$\mathcal{W}_l^{\rm ab}$};
   \draw (4.5,1.25) arc (270:90:.0625);
   \draw[->] (4.5,0) -- (5.5,0);
   \draw[->] (6,1) -- (6,.25);
   %\draw (6,.625) node [right] {$\alpha$};
   \draw (4,0) node {$\mathcal{W}_k^{\rm ab}$};
   \draw[<-] (4,.25) -- (4,1);
   %\draw (4,.625) node [left] {$\alpha$};
   \draw (4,1.25) node {$\mathcal{C}_k$};
   \draw[->] (4.5,1.25) -- (5.5,1.25);
   \draw (5,0) node [below] {${\rm Ver}_{E/F}$};
\end{tikzpicture}.
\end{center}
The reciprocity map sends one dimensional Artin $L$-functions to Tate $L$-functions for Gr\"o\ss encharakters of $\mathbb{A}_k^\times = {\rm GL}_1(\mathbb{A}_k)$. For every $v \in {\rm Val}_\mathfrak{p}(k)$, we retrieve the isomorphism
\begin{equation*}
   {\rm Art}_v : {\rm GL}_1(k_v) \xrightarrow{\sim}
    \mathcal{W}_{k_v}^{\rm ab},
\end{equation*}
of local class field theory.
}

\section{Reductive groups and representation theory}\label{rgrt}

Reductive groups are at the base of the theory, and were studied in detail by the Chevalley seminar and the S\'eminaire de G\'eom\'etrie Alg\'ebrique  (SGA 3). The scheme theoretic point of view gives a functor $\bf G$ that allows us to change the field or ring of definition. For practical purposes, as in the theory of Lie groups, becoming familiar with a list of examples works best for familiarizing oneself with the representation theory of reductive groups. For the reader unfamiliar with reductive groups, the examples of \S~\ref{examples} suffice on a first reading; for the basic notions of reductive groups we refer to \cite{SpBook} and  \cite{MiBook}.
  
We let $\bf G$ be a quasi-split connected reductive group defined over a global field $k$. We take $G = {\bf G}(\mathbb{A}_k)$ and $\Gamma = {\bf Z}_G(\mathbb{A}_k) {\bf G}(k)$. Langlands made a deep study of Eisenstein series in \cite{La1976}, providing a decomposition
\begin{equation*}
   \mathscr{L}^2(\Gamma \backslash G) =
   \mathscr{L}^2(\Gamma \backslash G)_{\rm disc} \oplus 
   \mathscr{L}^2(\Gamma \backslash G)_{\rm cont},
\end{equation*}
where we have the notion of cusps that live in the discrete spectrum and the process of parabolic induction which helps us understand the continuous spectrum. Langlands' study was initially over the field $\mathbb{R}$ and he later reformulated it in ad\`elic language for number fields. Over function fields, Harder made an initial study for split $\bf G$ in \cite{Ha1974}, which was completed by Morris in \cite{Mo1982}.

Let ${\bf P} = {\bf M}{\bf N}$ be a parabolic subgroup of $\bf G$ with Levi $\bf M$ and unipotent radical $\bf N$. We observe that $\bf M$ is a quasi-split connected reductive group over $k$ and we look at the space of functions
\begin{equation*}
   \mathscr{L}^2({\bf Z_M}(\mathbb{A}_k) {\bf M}(k) 
   									\backslash {\bf M}(\mathbb{A}_k)).
\end{equation*}
We say that an automorphic representation $\pi$ of ${\bf M}(\mathbb{A}_k)$ is cuspidal if for every $f$ in the space of $\pi$ and every parabolic subgroup ${\bf P}_0 = {\bf M}_0 {\bf N}_0$ of $\bf M$ we have that
\begin{equation*}
   \int_{{\bf N}_0(k) \backslash {\bf N}_0(\mathbb{A}_k)} f(ng) \, dn =0.
\end{equation*}

Langlands observed in \cite{LaCorvallisI} that every cuspidal automorphic representation $\pi$ of ${\bf M}(\mathbb{A}_k)$ arises as an irreducible constituent of the right regular representation of $\mathscr{L}^2(\Gamma \backslash M)$. Furthermore, every automorphic representation $\pi$ of ${\bf G}(\mathbb{A}_k)$ arises as a constituent of
\begin{equation*}
   {\rm Ind}_{{\bf P}(\mathbb{A}_k)}^{{\bf G}(\mathbb{A}_k)}(\rho),
\end{equation*}
with $\rho$ a cuspidal representation of ${\bf M}(\mathbb{A}_k)$, extended to ${\bf P}(\mathbb{A}_k)$ by taking it to be trivial on ${\bf N}(\mathbb{A}_k)$. Here, ${\rm Ind}$ denotes normalized parabolic induction, to which we will return locally in \S~\ref{admissible} and globally in \S~\ref{Ind:auto}. For now, it suffices for us to know that cuspidals are the building blocks of the theory and all automorphic representations $\pi$ are obtained by this process.

\subsection{Examples}\label{examples}
The group ${\rm GL}_1$ arises as the maximal torus (and Levi subgroup) $\bf T$ consisiting of diagonal matrices of the ambient group ${\bf G}={\rm SL}_2$. Automorphic representions of ${\rm GL}_1(\mathbb{A}_k)$ are Gr\"o\ss encharakters 
\begin{equation*}
   \chi : k^\times \backslash \mathbb{A}_k^\times \rightarrow \mathbb{C}^\times.
\end{equation*}
As observed by Langlands in \cite{La1967}, Eisenstein series in this setting allow us to study $L$-functions in concordance with the abelian case of Tate's thesis.

Let $n$ be a positive integer. We have in mind general linear groups and split classical groups when working with reductive groups in this article, namely:
\begin{equation*}
   \left.
      \begin{array}{cl}
         {\rm (a)} & {\bf G}_n = {\rm GL}_n \\
         {\rm (b)} & {\bf G}_n = {\rm SO}_{2n+1} \\
         {\rm (c)} & {\bf G}_n = {\rm Sp}_{2n} \\
         {\rm (d)} & {\bf G}_n= {\rm SO}_{2n}
      \end{array}
   \right\}
   \ 
   \begin{array}{c}
         \text{with their Levi subgroups} \\
         {\bf M} \cong {\rm GL}_{n_b} \times \cdots \times {\rm GL}_{n_1} \times {\bf G}_{n_0}
   \end{array}
\end{equation*}
$n = n_k + \cdots + n_1 + n_0$ and ${\bf G}_{n_0}$ a group of the same kind as ${\bf G}_{n_0}$. These reductive groups provide us with intuition into the general theory. And, indeed, their Lie algebras are of type $A_{n-1}$, $B_n$, $C_n$ and $D_n$, respectively. For the classical groups, i.e. cases (b)--(d), we have Levi subgroups with $n_0 = 0$ and we interpret ${\bf G}_0$ as empty.

In these examples, it is not hard to describe all parabolics ${\bf P} = {\bf M}{\bf N}$ explicitly over any ring $A$\footnote{Assume $A$ to be commutative with identity.}. For this, let $\mathfrak{M}_{s,t}(A)$ be the set of $s \times t$ matrices with entries in $A$ and let $I_s \in \mathfrak{M}_{s,s}(A)$ denote the identity matrix.

In case (a), we can precisely write the rational points of $\bf M$ and $\bf N$ as
\begin{align*}
   M & = \left\{ \left(X_{ij} \right) \vert \ X_{ij}\in \mathfrak{M}_{n_i,n_j}(A), \text{ with } X_{ii} \in {\rm GL}_{n_i}(A) \text{ and } X_{ij}=0 \text{ for } i \neq j \right\}, \\
   N & = \left\{ \left(X_{ij} \right) \vert \ X_{ij} \in \mathfrak{M}_{n_i,n_j}(A), X_{ii}=I_{n_i} \text{ and } X_{ij}=0 \text{ for } i>j \right\}.
\end{align*}
Parabolic subgroups $\bf P$ of ${\bf G}_n$, with $\bf M$ and $\bf N$ satisfying the above explicit description for every ring $A$, are the standard parabolics. From the theory of reductive groups, every parabolic subgroup is conjugate to a standard one.

Now assume we are in the case of ${\bf G}_n$ a classical group. Let $J_n \in \mathfrak{M}_{n,n}(A)$ be the matrix with $i,j$ entries given by $\delta_{i,n+1-j}$. If ${\bf G}_n$ is symplectic, i.e. ${\bf G}_n = {\rm Sp}_{2n}$, let
\begin{equation*}
   \Phi_n = \left(
      \begin{array}{cl}
         \ \ 0 & J_n \\
         -J_n & 0
      \end{array}
   \right).
\end{equation*}
And we write $\Phi_n = J_{n'}$ ,where $n' = 2n+1$ or $2n$ for odd or even orthogonal groups, respectively (for symplectic groups, set $n'=2n$). We can simplify things if we assume $-1 \neq 1$. Let us do so, although the theory is available without this restriction. Then the rational points, $G_n = {\bf G}_n(A)$, of a classical group are given by
\begin{equation*}
   G_n = \left\{ g \in {\rm SL}_{n'}(A) \vert \ {}^tg \Phi_n g = \Phi_n \right\}.
\end{equation*}

The standard parabolic subgroups ${\bf P}={\bf M}{\bf N}$ of ${\bf G}_n$ have rational points given as follows: for $M = {\bf M}(A)$ by the following set of block matrices
\begin{align*}
   M =
   \{  \left( X_{ij} \right) \vert \ &X_{b+1,b+1} \in G_{n_0}, X_{ij}=0 \text{ for } i \neq j, \text{ and } X_{ii}\in {\rm GL}_{n_i}(A), \\
   & X_{b'+1-i,b'+1-i} = J_{n_i} {}^tX_{i,i}^{-1} J_{n_i} \text{ for } 1 \leq i \leq b \},
\end{align*}
where $b' = 2b +1$ in case (b), and $b' = 2b$ in cases (c) and (d); the group of rational points $N$ of the unipotent radical $\bf N$ then consists of the corresponding block upper triangular matrices with block identity matrices along the diagonal.

\subsection{Orthogonal groups in characteristic $2$} While the symplectic group ${\rm Sp}_{2n}$ can be defined as before regardless of characteristic, the orthogonal groups require quadratic forms; especially ${\rm SO}_{2n}$. The previous simplification is possible because quadratic forms are equivalent to bilinear forms if $-1 \neq 1$.

For quadratic forms defined over a ring $A$, we refer to \cite{Kn1991}. In particular, let
\begin{align*}
   q_{2n+1}(x) = \sum_{i=1}^n x_i x_{2n+2-i} + x_{n+1}^2 \text{ and }
   q_{2n}(x) = \sum_{i=1}^n x_i x_{2n+1-i},
\end{align*}
for $x \in A^{n'}$, where $n' = 2n+1$ or $2n'$, respectively. Then ${\rm O}_{n'}$ is the orthogonal group of rank $n$ for the corresponding quadratic form $q_{n'}$.

For the special orthogonal groups, let $\mathbb{Z}_2(A)$ be the additive group of continuous maps ${\rm Spec}(A) \rightarrow \mathbb{Z}/2\mathbb{Z}$. There is the Dickson invariant, see [\emph{loc.\,cit.}], for the corresponding quadratic space:
\begin{equation*}
   {\rm O}_{n'}(A) \xrightarrow{D_{q_{n'}}} \mathbb{Z}_2(A) \longrightarrow 1,
\end{equation*}
Then, ${\rm SO}_{n'}(A)$ is defined to be the kernel of the epimorphism $D_{q_{n'}}$.

\subsection{Admissible representations of $\mathfrak{p}$-adic groups}\label{admissible}
The first article on $\mathfrak{p}$-adic groups was written by Mautner \cite{Ma1958}, who makes a study of spherical functions of ${\rm GL}_2$ and ${\rm PGL}_2$. This active area of current research has come a long way since then, and a modern introductory approach in accordance to our exposition can be found in \cite{BuHe2006}. The reader can assume that $\bf G$ is given by one of the examples of \ref{examples}. However, our framework is in the context of a general connected quasi-split reductive group $\bf G$, which means it has a Borel subgroup.

We fix a Borel subgroup ${\bf B}={\bf T}{\bf U}$ of $\bf G$. It is a minimal parabolic subgroup, with maximal torus $\bf T$ and unipotent radical $\bf U$. In our examples, $\bf B$ consists of upper triangular matrices and $\bf T$ consists of diagonal matrices. We assume all parabolic subgroups ${\bf P}={\bf M}{\bf N}$ are standard, by requiring ${\bf P} \supset {\bf B}$.

Given an algebraic group $\bf H$ defined over a non-archimedean local field $F$, we write $H={\bf H}(F)$ for its group of rational points. Then $G$ is a totally disconnected group, or, in other words, it is locally profinite. There exists then a basis of compact open subgroups (c.o.s.g. for short) $K$ of $G$.

We study representations $(\pi,V)$ of $G$, where $V$ is a complex (often infinite dimensional) vector space. A representation $(\pi,V)$ is smooth if for every $v \in V$, there is a c.o.s.g. $K$ of $G$ such that $\pi(K)v = v$. A smooth representation $(\pi,V)$ is admissible if for every c.o.s.g. $K$ the space of fixed vectors $V^K$ has
\begin{equation*}
   \dim_{\mathbb{C}} V^K < \infty.
\end{equation*}
All representations in the local Langlands conjecture are admissible, and this assumption allows us to define the trace of $\pi$.

Given a parabolic subgroup ${\bf P}={\bf M}{\bf N}$ of $\bf G$, the group $\bf M$ is itself quasi-split reductive. From the point of view of Haar measures, $G$ and $M$ are unimodular. However, this is no longer the case for $P$, which has non-trivial modulus character $\delta_P$ that allows us to pass between right Haar measures (which we use in all of our integrals) to left Haar measures. Let $\rho$ be an irreducible admissible representation of $M$, which we extend to one of $P$ by tensoring with the trivial representation of $N$. Then, parabolic induction defines a representation on $G$ via the right regular action on the space
\begin{align*}
   {\rm Ind}_P^G(\rho) = \left\{ f:G \rightarrow V \right. &\vert
   f(mng) = \delta_P^{1/2}(m)\rho(m)f(g), \, \forall m \in M, n \in N, \\
   &\left. \exists \text{ c.o.s.g. } K \text{ of } G \text{ such that } f(gk) = f(g), \, \forall k \in K \right\}.
\end{align*}
When $P$ and $G$ are clear from context, we often simply write ${\rm Ind}(\rho)$ and further abuse notation by identifying the induced representation with its space ${\rm Ind}(\rho)$.

Because we normalize by the character $\delta_P^{1/2}$, this parabolic induction is often called unitary parabolic induction. If $\rho$ is admissible, resp. unitary, then ${\rm Ind}(\rho)$ is also admissible, resp. unitary. Without the normalization, parabolic induction does not preserve unitarity. Furthermore, ${\rm Ind}(\rho)$ is of finite length, see \S~6 of \cite{CaNotes}. We also have the operation of taking the smooth dual, i.e., the contragredient representation $(\tilde{\rho},\tilde{V})$, which preserves admissibility and is respected by ${\rm Ind}$.

The building blocks of the theory are supercuspidal representations. These are irreducible admissible representations $\pi$ which cannot be obtained as a subquotient of a parabolically induced representation
\begin{equation*}
   %\pi \hookrightarrow
   {\rm Ind}_P^G(\rho),
\end{equation*}
with $\rho$ an admissible representation of $M$ and $M \neq G$. In our examples, supercuspidals can be constructed using the theory of types in \cite{BuKu1993,St2008}.

Let $(\pi,V)$ be an irreducible admissible representation of $G$. A matrix coefficient of $\pi$ is a function defined by
\begin{equation*}
   c_{v,\tilde{w}}(g) = \left\langle \pi(g) v, \tilde{w} \right\rangle, \ g \in G,
\end{equation*} 
where $v \in V$ and $\tilde{w} \in \tilde{V}$. With matrix coefficients we can define several important classes of representations in the theory.

An irreducible admissible representation $\pi$ of $G$ is supercuspidal if and only if all of its matrix coefficients have compact support. Discrete series representations are those whose matrix coefficients are in $\mathscr{L}^2(G)$. Tempered representations can be characterized by having their matrix coefficients in $\mathscr{L}^{2+\epsilon}(G)$ for every $\epsilon>0$.

Langlands classification states that every irreducible admissible representation $\pi$ of $G$ arises as a Langlands quotient of an induced representation
\begin{equation*}
   {\rm Ind}(\rho \cdot \chi),
\end{equation*}
where $\rho$ is a tempered representation of $M$ and $\chi$ is a character in the Langlands situation \cite{BoWa,Si1978}. What we have is the following diagram
\begin{equation*}
   \left\{ \begin{array}{c} \text{supercuspidal} \\
  				      \text{representations}
					\end{array} \right\}
   \subset
   \left\{ \begin{array}{c} \text{discrete} \\
   				     \text{series}
					\end{array} \right\}
	\subset
	\left\{ \begin{array}{c} \text{tempered} \\
  				      \text{representations}
					\end{array} \right\}
   \subset
	\left\{ \begin{array}{c} \text{admissible} \\
  				      \text{representations}
					\end{array} \right\}
\end{equation*}

\subsection*{Example}
Let $\delta$ be a unitary discrete series representation of ${\rm GL}_n(E)$. Zelevinsky classification \cite{Ze1980} tells us that it is a constituent of the induced representation
\begin{equation}\label{zclass}
   {\rm Ind} (\rho\nu^{-\frac{t-1}{2}} \otimes \cdots \otimes \rho\nu^{\frac{t-1}{2}}),
\end{equation}
where $\rho$ is a supercuspidal representation of ${\rm GL}_m(E)$, $m \vert n$, and $t$ is a positive integer. The representation $\delta$ is precisely the generic constituent, a notion towards which we now turn, and continue with further examples in \S~\ref{generic:examples}.

\subsection{$L$-groups and Satake parametrization}\label{Lgroups:Satake}
Given a quasi-split connected reductive group $\bf G$, Chevalley's theorem gives the existence of the dual group $\bf \widehat{G}$, and we take $\widehat{G}={\bf \widehat{G}}(\mathbb{C})$. Given a local or a global field $k$, we have the Langlands dual group
\begin{equation*}
   {}^LG = \widehat{G} \rtimes \mathcal{W}_k',
\end{equation*}
where $\mathcal{W}_k' = \mathcal{W}_k \times {\rm SL}_2(\mathbb{C})$ is the Weil-Deligne group.

Let now $F$ be a non-archimedean local field with ring of integers $\mathcal{O}_F$ and residue field $\mathbb{F} = \mathcal{O}_F / \mathfrak{p}_F$, with cardinality $q_F$. Since we assume $\bf G$ to be quasi-split over $F$, there is a Borel subgroup $\bf B$ of $\bf G$. If the group is split, then the groups are also defined over $\mathcal{O}_F$. If $\bf G$ is non-split quasi-split, one can always choose a connected group scheme ${\bf G}_0$ over $\mathcal{O}_F$, such that it has the same generic fiber and connected special fiber \cite{Ti1979}. Then, we can define the Iwahori subgroup as the pre-image in ${\bf G}_0(\mathcal{O}_F)$ of the $\mathbb{F}$-points of $\bf B$. Namely
\begin{equation*}
   \mathcal{I}_F = {\rm Ker}\left\{ {\bf G}_0(\mathcal{O}_F) \rightarrow {\bf B}(\mathbb{F}) \right\}.
\end{equation*}

Let
\begin{equation*}
   \mathfrak{Rep}_{\rm nr}(G) = \left\{ \begin{array}{c}
   	\text{isomorphism classes of smooth irreducible } \\
   	\text{ representations } \pi \text{ of } G \text{ such that }
   	\pi^{\mathcal{I}_F} \neq 0
				\end{array} \right\}.
\end{equation*}
It is well known that this class of representations satisfies the local Langlands correspondence \cite{Sa1963}, a result that is known as the Satake parametrization:
\begin{equation*}
   \pi \in \mathfrak{Rep}_{\rm nr}(G) \quad \longleftrightarrow \quad
   \phi_\pi: \mathcal{W}_F' \rightarrow {}^LG,
\end{equation*}
where $\phi_\pi$ is the local Langlands parameter attached to $\pi$. It is represented by
\begin{equation*}
   \phi_\pi({\rm Frob}_F) = s(\pi),
\end{equation*}
a Satake parameter $s(\pi)$ inside the variety of semi-simple elements $[\widehat{G}^{\mathcal{I}_F} \rtimes {\rm Frob}_F]_{\rm ss}$ modulo conjugation under $\widehat{G}^{\mathcal{I}_F}$. A general construction can be found in \cite{HaSatake}. 

\begin{definition}\label{nr:L}
Let $r$ be an analytic finite dimensional representation of ${}^LG$ and let $\pi \in \mathfrak{Rep}_{\rm nr}(G)$. The corresponding local $L$-function is defined by
\begin{equation*}
   L(s,\pi,r) = \dfrac{1}{\det(I - r(s(\pi))q_F^{-s})}, \ s \in \mathbb{C}.
\end{equation*}
\end{definition}

We also note the fact that unramified representations are obtained via parabolic induction from an unramified character of the maximal torus. This is a basic result of Borel and of Casselman \cite{Ca1979}.

\begin{theorem}
Let $\bf G$ be a quasi-split connected reductive group over $F$ with Borel subgroup ${\bf B} = {\bf T}{\bf U}$. If $\pi \in \mathfrak{Rep}_{\rm nr}(G)$, then there exists an unramified character $\chi$ of $T$ such that $\pi$ is the unique irreducible generic constituent of
\begin{equation*}
   {\rm Ind}(\chi)
\end{equation*}
\end{theorem}

\subsection*{Examples} Our examples of type (a)--(d) are all split. Because of this, given our rank $n$ group ${\bf G}_n$, we only require the connected component of the Langlands dual group in the parametrization:
\begin{equation*}
   {}^LG_n^\circ = \widehat{G}_n =
   \left\{
      \begin{array}{ll}
         {\rm GL}_n(\mathbb{C}) & 
         \text{ if  } {\bf G}_n = {\rm GL}_n \\
         {\rm Sp}_{2n} (\mathbb{C}) & 
         \text{ if  }{\bf G}_n = {\rm SO}_{2n+1} \\
         {\rm SO}_{2n+1}(\mathbb{C}) & 
         \text{ if  } {\bf G}_n = {\rm Sp}_{2n} \\
         {\rm SO}_{2n}(\mathbb{C}) & 
         \text{ if  }{\bf G}_n= {\rm SO}_{2n}
      \end{array}
   \right.
   .
\end{equation*}
We observe that the general linear group ${\rm GL}_n$ and the even special orthogonal group ${\rm SO}_{2n}$ are self dual.

Fix a classical group ${\bf G}_n$ and let $\pi \in \mathfrak{Rep}_{\rm nr}(G_n)$. Then, there exist unramified characters $\chi_1, \ldots, \chi_n$ of ${\rm GL}_1(F)$ such that
\begin{equation*}
   \pi \hookrightarrow {\rm Ind}(\chi_1 \otimes \cdots \otimes \chi_n)
\end{equation*}
is the unique irreducible generic constituent. The Satake parameters for ${\bf G}_n = {\rm GL}_n$ are then given by the diagonal (hence semisimple) matrix
\begin{equation*}
   s(\pi) = 
   {\rm diag}(\chi_1(\varpi_F), \ldots, \chi_n(\varpi_F)) 
   \in {\rm GL}_n(\mathbb{C}).
\end{equation*}
If ${\bf G}_n = {\rm SO}_{2n+1}$ or ${\rm SO}_{2n}$, then
\begin{equation*}
   s(\pi) = 
   {\rm diag}(\chi_1(\varpi_F), \ldots, \chi_n(\varpi_F),
   \chi_n^{-1}(\varpi_F), \ldots, \chi_1^{-1}(\varpi_F)) 
   \in {}^LG_n^\circ.
\end{equation*}
And if ${\bf G}_n = {\rm Sp}_{2n}$, then
\begin{equation*}
   s(\pi) = 
   {\rm diag}(\chi_1(\varpi_F), \ldots, \chi_n(\varpi_F), 1,
   \chi_n^{-1}(\varpi_F), \ldots, \chi_1^{-1}(\varpi_F)) 
   \in {}^LG_n^\circ.
\end{equation*}

\subsection{Induction and automorphic representations}\label{Ind:auto}
We use the guidance of Langlands, who observed that the process of parabolic induction allows us to obtain all automorphic representations from cuspidal representations \cite{LaCorvallisI}.

We are in the situation of a quasi-split connected reductive group $\bf G$, defined over a global field $k$. And, ${\bf P} = {\bf M}{\bf N}$ is a parabolic subgroup. We write $G_v$ for ${\bf G}(k_v)$ and similarly for $P_v$, $M_v$ and $N_v$. For each $v \in {\rm Val}_\mathfrak{p}(k)$, we have the ring of integers $\mathcal{O}_v$ of $k_v$. The group $K_v = {\rm G}(\mathcal{O}_v)$ is defined for almost every $v \in {\rm Val}_\mathfrak{p}(k)$, and it is a hyperspecial maximal compact open subgroup of $G_v$.

Let $\rho$ be a cuspidal automorphic representation of ${\bf M}(\mathbb{A}_k)$ extended to ${\bf P}(\mathbb{A}_k)$ by taking it to be trivial on ${\bf N}(\mathbb{A}_k)$. Every automorphic representation can be written as a restricted tensor product
\[ \rho = \otimes' \rho_v, \]
such that $\rho_v \in \mathfrak{Rep}_{\rm nr}(G_v)$ for almost all $v \in {\rm Val}(k)$, see \cite{Fl1979}. The decomposition is made by fixing a family of vectors $x_v$ in the space of $\rho_v$, fixed by $K_v = {\bf M}(\mathcal{O}_v)$, for almost all $v$. Then we can form the restricted tensor product
\[ {\rm Ind}_{{\bf P}({\mathbb{A}_k})}^{{\bf G}({\mathbb{A}_k})} \rho
   = \otimes' {\rm Ind}_{P_v}^{G_v} \rho_v. \]
The restricted tensor product is taken with respecto to a family of nonzero vectors $f_v^0 $ in the space of ${\rm Ind}_{P_v}^{G_v} \rho_v$. We fix such a family by taking $S$ to be a finite set of places of $k$, $S \supset {\rm Valp}_\infty(k)$, such that $\rho_v$ is unramified for $v \notin S$; we can further assume that $K_v = {\bf G}(\mathcal{O}_v)$. We have a unique constituent
\[ \pi_v^\circ \text{ of } {\rm Ind}_{P_v}^{G_v} \rho_v, \ v \notin S, \]
containing the trivial representation of $K_v$. Then, for each $v \notin S$, the function $f_v^0$ is taken in the space of $\pi_v^\circ$ to be such that $f_v^0(y_v) = x_v$ for all $y_v \in K_v$. Then every function in $V_\pi$ is of the form
\[ f = \otimes'_v f_v,  \quad f_v = f_v^0, \ v \notin S', \]
where $S'$ is a set containing $S$.

The main result of \cite{LaCorvallisI}, is that a representation $\pi$ of ${\bf G}(\mathbb{A}_k)$ is an automorphic representation if and only if $\pi$ is a constituent of
\[ {\rm Ind}_{{\bf P}({\mathbb{A}_k})}^{{\bf G}({\mathbb{A}_k})} \rho \]
for a parabolic ${\bf P} = {\bf M}{\bf N}$ of $\bf G$ and a cuspidal representation $\rho$ of ${\bf M}(\mathbb{A}_k)$. The representation can be expressed as a tensor product
\[ \pi = \otimes' \pi_v, \]
where $\pi_v = \pi_v^\circ$ for almost all $v$.
 
\subsection{Generic representations}
 Fixing a Borel subgroup $\bf B$ of $\bf G$ gives a pinning of the roots, and in particular determines the set of positive roots $\Phi^+$. We use Bourbaki notation, chapter VI of \cite{Bou}, when working with root systems. For every $\alpha \in \Phi^+$ there is a rank one group ${\bf G}_\alpha$ obtained via restriction of scalars from either ${\rm SL}_2$ or ${\rm SU}_3$. In the former case we have that the unipotent subgroup $N_\alpha = U_\alpha$ is a one parameter additive subgroup and in the latter case we have $N_\alpha = U_\alpha U_{2\alpha}$. Our examples are of split groups, in which the latter possibility does not occur.

Let $\Delta$ denote the simple roots. We have the surjective morphism
\begin{equation}\label{psiU}
   U \twoheadrightarrow U / \prod_{\alpha \in \Phi^+ \setminus \Delta} N_\alpha \simeq \prod_{\alpha \in \Delta} U_\alpha,
\end{equation}
Given a non-trivial additive character $\psi: A \rightarrow \mathbb{C}^\times$, on a topological ring $A$, we extend it to one of $U = {\bf U}(A)$ via the morphism \eqref{psiU}
\begin{equation*}
   \psi(u) = \prod_{\alpha \in \Delta} \psi(u_\alpha).
\end{equation*}
In practice, we know how to pass from an arbitrary character $\chi$ of $U$ to one of the form just constructed.

\subsection*{Example} For ${\rm GL}_n$, we fix the Borel subgroup ${\bf B} = {\bf T}{\bf U}$ of upper triangular matrices. The maximal torus $\bf T$ consists of diagonal matrices and the unipotent radical $\bf U$ has rational points given by
\begin{equation*}
   U = {\bf U}(A) = \left\{ u = \left(x_{ij} \right) \vert \ x_{ij} \in A, x_{ii}=1 \text{ and } x_{ij}=0 \text{ if } i<j \right\} \subset {\rm GL}_n(A). 
\end{equation*}
From the non-trivial character $\psi:F \rightarrow \mathbb{C}^\times$, we extend $\psi$ to $u \in {\bf U}(A)$ by setting
\begin{equation*}
   \psi(u) = \prod_{i=1}^{n-1} \psi(x_{i,i+1}).
\end{equation*}
It is a good excercise to write the corresponding explicit matrix realization for the classical groups.

Let us now return to the notion of genericity, first, for local representations. Fix a non-archimedean local $F$ and a non-trivial character $\psi: F \rightarrow \mathbb{C}^\times$. We say that an irreducible admissible representation $(\pi,V)$ of ${\bf G}(F)$ is generic if there exists a continuous functional
\begin{equation*}
   \lambda: V \rightarrow \mathbb{C},
\end{equation*}
such that
\begin{equation*}
   \lambda(\pi(u)v) = \psi(u) \lambda(v), \text{ for } u \in U.
\end{equation*}
Such a $\lambda$, when it exists, is called a Whittaker functional. We have the follwing multiplicity one theorem due to Shalika \cite{ShJ1974}.

\begin{theorem}
   Let $(\pi,V)$ be a generic representation of $G = {\bf G}(F)$. Let $W$ denote the space of all Whittaker functionals, then
\begin{equation*}
   \dim_\mathbb{C} W = 1.
\end{equation*}
\end{theorem}

Now, for a global field $k$. Let $\pi = \otimes' \pi_v $ be a cuspidal automorphic representation of ${\bf G}(\mathbb{A}_k)$ and $\psi: k \backslash \mathbb{A}_k \rightarrow \mathbb{C}^\times$ a non-trivial character. We say that $\pi$ is globally generic if there is a cusp form $\varphi$ in the space of $\pi$ such that
\begin{equation*}
   W_{\varphi}(g) = \int_{{\bf U}(k) \backslash {\bf U}(\mathbb{A}_k)} \varphi(ug) \overline{\psi}(u) \, du \neq 0.
\end{equation*}
A globally generic representation $\pi \otimes' \pi_v$, is locally generic. Indeed, if we write the cusp form as a restricted product $\varphi = \otimes' \varphi_v$, ranging over all places $v$ of $k$, then we have
\begin{equation*}
   W_\varphi(g) = \prod_v W_{\varphi_v}(g_v) = \prod_v \lambda_v(\pi_v(g_v)\varphi_v),
\end{equation*}
where each $\lambda_v$ is a Whittaker functional for $\pi_v$.

\subsection{Unitary generic representations of the classical groups}\label{generic:examples}
Let ${\bf G}_n$ be either a general linear group or a classical group defined over a non-archimedean local field $F$. We have the following result concerning the unitary spectrum.

\begin{theorem}\label{lmt:genericU}
Let $\pi$ be an irreducible generic unitary admissible representation of ${\bf G}_n(F)$, then $\pi$ arises as arises as the irreducible quotient of a parabolically induced representation
\begin{equation}\label{Langlands:class:classical}
   {\rm Ind}(\tau_{b} \nu^{r_b} \otimes \cdots \otimes \tau_{1} \nu^{r_1} \otimes \pi_0).
\end{equation}
The induction is with respect to a parabolic subgroup with Levi ${\bf M} \cong {\rm GL}_{n_b} \times \cdots \times {\rm GL}_{n_1} \times {\bf G}_{n_0}$. For each $1 \leq i \leq b$, $\tau_i$ is a tempered unitary representation of ${\rm GL}_{n_i}(F)$, and the Langlands parameters satisfy $r_b \geq \cdots \geq r_1 \geq 0$. In the case of a general linear group, the representation $\pi_0$ of ${\bf G}_{n_0}(F)$ is quasi-tempered of the form $\pi_0' \nu^{r_0}$ with $\pi_0'$ tempered and $r_1 \geq r_0 \geq 0$. In the case of a classical group, $\pi_0$ is tempered generic. Furthermore
\begin{equation*}
   r_b < \left\{ \begin{array}{cl}
   						1 & \text{ if } {\bf G}_l = {\rm SO}_{2n} \\
						 	1/2 & \text{ otherwise }
				\end{array} . \right.
\end{equation*}
\end{theorem}

The bound $r_b < 1/2$ for ${\rm GL}_n$ is known as the trivial Ramanujan bound, obtained when viewed from a global perspective \cite{JaSh1981}. For the classical groups, a description of the generic unitary spectrum in characteristic zero is given in \cite{LaMuTa2004}. This result can be transferred to characteristic $p$ using the technique of Ganapathy on the Kazhdan transfer for close local fields \cite{Ga2015}, and we do so in \cite{LoRationality}.

\subsection{Intertwining operators and Plancherel measures}\label{inter}
Harish-Chandra studied the Plancherel theorem over the reals as well as the $\mathfrak{p}$-adics. We look at $\mu$-factors which make an appearence in the local Langlands conjecture (see \S~\ref{loc:L:classical}), in addition to their original role in the Plancherel theorem \cite{Wa2003}. 

Fix a pair of quasi-split reuductive groups $({\bf G},{\bf M})$ defined over a non-archimedean local field $F$, with $\bf M$ a maximal Levi subgroup of $\bf G$. We let
\begin{equation} \label{w_0}
   w_0 = w_l w_{l,\theta},
\end{equation}
where $w_l$ and $w_{l,\theta}$ are the longest elements in $G$ and in $M$, respectively. We have that conjugation by $w_0$ defines a Levi subgroup ${\bf M}_{w_0}$. Given a representation $\xi$ on $M$, we let $w_0(\xi)$ denote the representation on $M_{w_0}$ defined by $w_0(\xi)(m') = \xi(w_0^{-1}m'w_0)$.

Let $X_{\rm nr}(\bf M)$ be the group of unramified characters $\chi:M \rightarrow \mathbb{C}$; it is a complex algebraic variety. Given an irreducible admissible represenation $(\pi,V)$ of $M$, and a $\chi \in X_{\rm nr}(\bf M)$ we have an induced representation
\begin{equation*}
   {\rm I}(\chi,\pi) = {\rm Ind}_P^G(\pi \cdot \chi).
\end{equation*}
Furthermore, we have an intertwining operator
\begin{equation*}
   {\rm A}(\chi,\pi,w_0) : {\rm I}(\chi,\pi) \rightarrow {\rm I}(w_0(\chi),w_0(\pi)),
\end{equation*}
which is defined via the principal value integral
\begin{equation*}
   {\rm A}(\chi,\pi,w_0)f(g) = \int_{N_{w_0}} f(w_0^{-1}ng) \, dn.
\end{equation*}
It extends to a rational operator on $\chi \in X_{\rm nr}({\bf M})$, see \S~IV of \cite{Wa2003}.

Plancherel measures are defined via the composite of two intertwining operators
\begin{equation*}
    \mu(\chi,\rho,\psi)^{-1} = 
    A(w_0(\chi),w_0(\pi),\psi) \circ A(\chi,\pi,\psi).
\end{equation*}
We obtain this way a a scalar-valued rational function $\mu(\chi,\pi,\psi)$ on the variable $\chi \in X_{\rm nr}({\bf M})$.

\subsection*{Examples} It is customary to write $\nu$ for the unramified character $\left| \det(\cdot) \right|$ of ${\rm GL}_n$, and we do so. We note that every unramified character $\chi \in X_{\rm nr}({\rm GL}_n)$ is of the form $\chi = \nu^s$, as $s$ ranges through $\mathbb{C}$. In case (a) of examples~\ref{examples}, we obtain functions of one complex variable by taking the representation $\rho = \tau \otimes \tilde{\pi}$ of $M = {\rm GL}_m(F) \times {\rm GL}_n(F)$ and setting
\begin{equation*}
   {\rm I}(s,\rho) = {\rm Ind}_P^G(\tau \nu^{s/2} \otimes \tilde{\pi}\nu^{-s/2}).
\end{equation*}
For the classical groups, i.e., cases (b)--(d), we take the representation $\rho = \tau \otimes \tilde{\pi}$ of $M = {\rm GL}_m(F) \times {\rm G}_n$ and set
\begin{equation*}
   {\rm I}(s,\rho) = {\rm Ind}_P^G(\tau \nu^{s} \otimes \tilde{\pi}).
\end{equation*}
With these normalizations, the intertwining operator is of the form
\begin{equation*}
   {\rm A}(s,\rho,w_0) : {\rm I}(s,\rho) \rightarrow {\rm I}(-s,\tilde{w}_0(\rho)).
\end{equation*}
While $M_{w_0} \neq M$ in general for ${\rm GL}_l$, it is the case that $M_{w_0} = M$ for the classical groups. We also observe that we have rational functions
\begin{equation*}
   \mu(s,\rho,\psi) \in \mathbb{C}(q^{-s}).
\end{equation*}

\subsection{Partial $L$-functions} Let $\bf G$ be a quasi-split connected reductive group defined over a global field $k$. It is a basic property that an automorphic representation $\pi$ of ${\bf G}(\mathbb{A}_k)$, has a decomposition
\begin{equation*}
 \pi = \otimes' \pi_v,
\end{equation*}
so that there is a finite set $S$ of places of $k$ (containing the archimedean places in the case of a number field), for which $\pi_v \in \mathfrak{Rep}_{\rm nr}({\bf G}(k_v))$, for $v \notin S$ \cite{Fl1979}.

\begin{definition}\label{partial:L}
Let $r$ be an analytic finite dimensional representation of ${}^LG$. Let $\pi = \otimes' \pi_v$ be an automorphic representation of ${\bf G}(\mathbb{A}_k)$ and $S$ a finite set of places of $k$ such that $\pi_v \in \mathfrak{Rep}_{\rm nr}({\bf G}(k_v))$, for all $v \notin S$. Then, partial $L$-functions are defined by
\begin{equation*}
   L^S(s,\pi,r) = \prod_{v \notin S} L(s,\pi_v,r_v).
\end{equation*}
\end{definition} 

An basic property of partial automorphic $L$-functions is that $L^S(s,\pi,r)$ converges for $\Re(s) \gg 0$, see \cite{Bo1979}. However, it is much more difficult to obtain a meromorphic continuation and functional equation. These last two properties, we address in many cases via the Langlands-Shahidi method, to which we next turn.

\section{The Langlands-Shahidi method}\label{LSmethod}

\subsection{}\label{LS:intro}
Fix a pair of quasi-split connected reductive groups $({\bf G},{\bf M})$, with $\bf M$ a maximal Levi subgroup of $\bf G$. Locally, the Langlands-Shahidi method applies to:
\begin{itemize}
   \item[(i)] a generic representation $\pi$ of $M$;
   \item[(ii)] a constituent $r$ of the adjoint action of ${}^LM$ on ${}^L\mathfrak{n}$, where $\mathfrak{n}$ is the Lie algebra of the unipotent radical of $\bf N$;
   \item[(iii)] a non-trivial additive character $\psi$.
\end{itemize}
 
Given a local triple $(\pi,r,\psi)$, with $M = {\bf M}(F)$, the Langlands-Shahidi machinery produces rational functions
\begin{equation*}
   \gamma(s,\pi,r,\psi) \in \mathbb{C}(q^{-s}),
\end{equation*}
of a complex variable $s \in \mathbb{C}$. Furthermore, if $\pi$ is tempered, then the $\gamma$-factor can be written as a product of a monomial $\varepsilon(s,\pi,r,\psi)$ and a ratio of two monic polynomias on $q^{-s}$
\begin{equation}\label{temp:3}
   \gamma(s,\pi,r,\psi) = \varepsilon(s,\pi,r,\psi) \dfrac{L(1-s,\tilde{\pi},r)}{L(s,\pi,r)}.
\end{equation}
This definig relationship is such that for $\pi \in \mathfrak{Rep}_{\rm nr}(M)$ it agrees with the definition given by the Satake correspondence, where we have $\varepsilon(s,\pi,r,\psi)=1$. It is extended from tempered generic to all generic $\pi$ via Langlands classification and a multiplicativity property, explained in \S~\ref{loc:LS}.

A global triple $(\pi,r,\psi)$, with $M = {\bf M}(\mathbb{A}_k)$, would be one with $\pi$ a globally generic cuspidal automorphic representation of $M$, $r$ an irreducible constituent of the adjoint action of ${}^LM$ on ${}^L\mathfrak{n}$, and $\psi = \otimes' \psi_v$ an additive character of $k \backslash \mathbb{A}_k$. We let $S$ be a finite set of places of $k$ such that $\pi_v$, $k_v$ and $\psi_v$ are unramified for $v \notin S$, were we have the partial $L$-function $L^S(s,\pi,r)$ defined in the previous section.

The Langlands-Shahidi machinery proves the functional equation for global triples $(\pi,r,\psi)$
\begin{equation}\label{gamma:fe}
   L^S(s,\pi,r) = \prod_{v \in S} \gamma(s,\pi_v,r_v,\psi_v) L^S(1-s,\tilde{\pi},r).
\end{equation}
Furthermore, with the local theory in place, we have completed $L$-functions and global root numbers
\begin{equation}\label{comp:L3}
   L(s,\pi,r) = \prod_v L(s,\pi_v,r_v) \text{ and } 
   \varepsilon(s,\pi,r) = \prod_v \varepsilon(s,\pi_v,r_v,\psi_v).
\end{equation}

\subsection*{Examples} Suppose we are in cases (a)--(d) of examples~\ref{examples}. The Langlands-Shahidi method leads us to automorphic $L$-functions for Rankin-Selberg products
\begin{equation*}
   L(s,\pi \times \tau),
\end{equation*}
where $\pi$ and $\tau$ are globally generic cuspidal automorphic representations of ${\bf G}_n(\mathbb{A})$ and ${\rm GL}_m(\mathbb{A}_k)$, respectively. The underlying analytic representation of ${}^LM$ is given by $r = \rho_n \otimes \rho_m$, where for every positive integer $l$, we let $\rho_l$ denote the standard representation of ${\rm GL}_l(\mathbb{C})$. This is one of two possible irreducible constituents of the adjoint action that arise in cases (b)--(d), the case (a) giving only Rankin-Selberg products.

We note that for ${\rm GL}_n$, all cuspidal automorphic representations are globally generic. This is no longer the case for the classical groups, as we mention in \S~\ref{GenRamanujan}, however, locally the tempered $L$-packet conjecture of Shahidi is known for the classical groups \cite{Ar2013,  GaVa2017}, which allows us to obtain all local $L$-functions from the generic case.

In cases (b)--(d) we also have the related $L$-functions $L(s,\tau,r)$ for a globally generic cuspidal automorphic representation $\tau$ of ${\rm GL}_m(\mathbb{A}_k)$, and we can have $r= \rho_m$, $r = {\rm Sym}^2 \rho_m$ or $\wedge^2 \rho_m$. These arise in the Siegel Levi case of ${\bf M} \cong {\rm GL}_m$ as a maximal Levi subgroup of a classical group ${\bf G}_m$ of rank $m$.

\subsection{Induced representations and genericity}

Fix a pair of quasi-split reuductive groups $({\bf G},{\bf M})$, with $\bf M$ a maximal Levi subgroup of $\bf G$. We consider the cases of $A=F$, a non-archimedean local field, or $A = \mathbb{A}_k$, the ring of ad\`eles of a global field. We write $G = {\bf G}(A)$, the rational points, and similarly for $M$.

Given a non-trivial characrter $\psi$ of $A$, we first use equation \eqref{psiU} to obtain a character $\psi_M$ of $U_M$, then we then require the following relation on a character $\psi$ of $U$:
\begin{equation*}
   \psi(u) = \psi_M(w_0^{-1}uw_0), \ u \in U_M = U \cap M,
\end{equation*}
where $w_0 = w_l w_{l,M}$ is the product of the longest Weyl group element $w_l$ of $G$ and the longest Weyl group element $w_{l,M}$ of $M$.

Let ${\bf P} = {\bf M}{\bf N}$ be a maximal parabolic subgroup of $\bf G$. The group ${\bf M} / {\bf Z}_{\bf G}$ has a one dimensional split center, hence ${\rm Hom}_F({\bf M}/{\bf Z}_{\bf G},{\rm GL}_1)$ is a rank one free $\mathbb{Z}$-module. Let
\begin{equation*}
   \delta = \det({\rm Ad}_{\bf M} \vert_\mathfrak{n}) \in 
   {\rm Hom}_F({\bf M}/{\bf Z}_{\bf G},{\rm GL}_1).
\end{equation*}
Now, in order to have $L$-functions of a complex variable, we let $\alpha$ be the simple root such that the maximal parabolic ${\bf P}$ corresponds to the subset of simple roots $\theta = \Delta - \left\{ \alpha \right\}$. Let $\delta_{\alpha}$ be an unramified character of $M$ given by
\begin{equation}\label{deltatilde}
   \delta_\alpha: M \xrightarrow{\delta^l} {\rm GL}_1(A) 
   \xrightarrow{\left| \cdot \right|_A} \mathbb{C}^\times, \ 
   l = \left\langle \delta, \alpha^\vee \right\rangle^{-1},
\end{equation}
where $\left\langle \cdot\, , \cdot \right\rangle$ is the perfect paring between roots and co-roots.

Let $\pi$ be a generic representation of $M$. We write
\begin{equation*}
   {\rm I}(s,\pi) = {\rm I}(\delta_{\!\alpha}^{s},\pi), \ s \in \mathbb{C},
\end{equation*}
with $\delta_\alpha$ as in equation~\eqref{deltatilde}; its corresponding space is denoted by ${\rm V}(s,\pi)$. %Furthermore, we write ${\rm I}(\pi)$ for ${\rm I}(0,\pi)$.

\subsection*{Examples} Suppose we are in the case of a maximal Levi subgroup
\[ {\bf M} \cong {\rm GL}_m \times {\bf G}_n \subset {\bf G}_l, \]
$l = m + n$, and both $m$, $n > 0$. If we are given $\psi$-generic representation $\tau$ of ${\rm GL}_m$ and $\pi$ of $G_n$, then we form the $\psi_M$-generic representation
\[ \rho = \tau \otimes \tilde{\pi} \text{ of } M. \]
The use of the contragredient representation $\tilde{\pi}$ is made in order to obtain Rankin-Selberg product $L$-functions from the local Langlands-Shahidi method. In all of these cases we have
\[ {\rm I}(s,\rho) = {\rm I}((\tau \cdot \nu^{s})\otimes \tilde{\pi}), \ s \in \mathbb{C}, \]
where $\nu$ denotes the unramified character by composing $\det$ with the absolute value $\left| \cdot \right|$ for ${\rm GL}_m$.

We have the following result, Theorem~2.1 of \cite{LoRationality}.

\begin{theorem}\label{Whittaker:Ind}
Let $(F,\pi,\psi)$ be a local triple in the Langlands-Shahidi setting. Then the induced representation ${\rm I}(s,\pi)$ is also generic. It has a Whittaker functional $\lambda_\psi(s,\pi)$ which is a polynomial on $\left\{ q^s, q^{-s} \right\}$. 
\end{theorem}

\subsection{Local Langlands-Shahidi method}\label{loc:LS}

Fix a pair of quasi-split reuductive groups $({\bf G},{\bf M})$, $\bf M$ a maximal Levi subgroup of $\bf G$, defined over a non-archimedean local field $F$. Given a local triple $(\pi,r,\psi)$, the Langlands-Shahidi method produces a local coefficient $C_\psi(s,\pi,w_0)$, which we combine with an induction step to produce $\gamma$-factors, $L$-functions and $\varepsilon$-factors.

From Theorem~\ref{Whittaker:Ind}, we have a Whittaker functional $\lambda_\psi(s,\pi)$ on ${\rm I}(s,\pi)$, and from \S~\ref{inter} we have an intertwining operator
\begin{equation*}
   {\rm A}(s,\pi,w_0) : {\rm I}(s,\pi) \rightarrow {\rm I}(w_0(\delta_{\!\alpha}^{s}),w_0(\pi)),
\end{equation*}
where in our examples, the unramified character $\delta_{\!\alpha}^{s}$ corresponds to $-s$. The crucial result is Shalika's multiplicity one theorem for Whittaker functionals \cite{ShJ1974}.

\begin{theorem}[Shalika]
An irreducible admissible representation $\Pi$ of $G$ is generic if and only if its space of Whittaker functionals is one dimensional.
\end{theorem}

We can then define the Langlands-Shahidi local coefficient via the formula
\[ \lambda_\psi(s,\pi) = C_\psi(s,\pi,w_0) \, \lambda_\psi(\delta_{\!\alpha}^{s},w_0(\pi)) \circ {\rm A}(s,\pi,w_0). \]
We have seen that the intertwining operator is rational and the Whittaker functional is a Laurent polynomial, we conclude that
\[ C_\psi(s,\pi,w_0) \in \mathbb{C}(q^{-s}). \]

We now observe that the adjoint action of ${}^LM$ on ${}^L\mathfrak{n}$, decomposes into a finite number $m(R)$ of irreducible components. In our examples, we actually have $m(R) = 1$ or $2$. Then, the following result, see for example Lemma~2.4 of \cite{LoLS}, tells us that we can inductively reduce to the case where the adjoint representation is irreducible.

\begin{lemma}\label{induction}
Let $({\bf G}, {\bf M})$ be a Langlands-Shahidi pair of quasi-split reductive groups, with ${\bf M}$ a maximal Levi subgroup of $\bf G$. Let $R = \oplus_{i=1}^{m(R)} r_i$ be the adjoint action of ${}^LM$ on ${}^L\mathfrak{n}$. For each $i > 1$, there exists a pair $({\bf G}_i,{\bf M}_i)$ such that the corresponding adjoint action of ${}^LM_i$ on ${}^L\mathfrak{n}_i$ decomposes as
\begin{equation*}
   R' = \bigoplus_{j=1}^{m(R')} r_j' \text{  with  } m(R') < m(R),  
\end{equation*}
and $r_{m(R)} = r_1'$.
\end{lemma}

We thus are able to obtain the individual $\gamma$-factors, which lead towards local $L$-functions and $\varepsilon$-factors, as observed in \S~\ref{LS:intro}. They satisfy a number of local properties, and we highlight two of them starting with their compatibility with Artin $\gamma$-factors, proved in \cite{KeSh1988} for non-archimedean local fields of characteristic $0$ and extended to characteristic $p$ in \cite{LoLS}.

\begin{theorem}
Let $(F,\pi,\psi)$ be a local triple such that $\pi \in \mathfrak{Rep}_{\rm nr}(M)$. Let 
\[ \phi_\pi : \mathcal{W}_F' \rightarrow {}^LM \]
be the Langlands parameter corresponding to $\pi$ via the Satake parametrization, see \S~\ref{Lgroups:Satake}. Then
   \begin{equation*}
      \gamma(s,\pi,r ,\psi) = \gamma(s,r \circ \phi_\pi,\psi).
   \end{equation*}
\end{theorem}

Local $\gamma$-factors satisfy a very interesting multiplicativity property, depending on cocycle relationships of Weyl group elements. For our examples, we can make this property explicit as follows.

\subsection*{Multiplicativity}
Let $({\bf G}, {\bf M})$ be a pair such that ${\bf G} = {\bf G}_l$ is a general linear group or a classical group of rank $l = m + n$, with $\bf M \cong {\rm GL}_m \times {\bf G}_n$ and both $m$, $n \geq 1$. Let $\tau$ be a generic representation of ${\rm GL}_m(F)$, such that it is a constituent of
\begin{equation}\label{gl:ind:rep}
   {\rm Ind}(\tau_1 \otimes \cdots \otimes \tau_e),
\end{equation}
where the $\tau_i$'s are generic representations of ${\rm GL}_{m_i}(F)$, $m = m_1 + \cdots m_e$. Let $\pi$ be a generic representation of ${\bf G}_n(F)$, such that it is a constituent of
\[ {\rm Ind}(\pi_1 \otimes \cdots \otimes \pi_f \otimes \pi_0), \]
where the $\pi_j$'s are generic representations of ${\rm GL}_{n_j}(F)$, $1 \leq j \leq f$, and $\pi_0$ is a generic representation of ${\bf G}_{n_0}(F)$, $n = n_1 + \cdots + n_f + n_0$, with ${\bf G}_{n_0}$ a classical group of the same type as ${\bf G}_{n}$.

The adjoint representation $R$ factors into two irreducible constituents $r_1$ and $r_2$. The first constituent is $r_1 = \rho_m \otimes \tilde{\rho}_n$, hence we obtain $\gamma(s,\tau \otimes \tilde{\pi},r_1,\psi) = \gamma(s,\tau \times \pi,\psi)$, a Rankin-Selberg $\gamma$-factor. If ${\bf G}_l = {\rm GL}_l$, then
\[ \gamma(s,\tau \times \pi,\psi)
   = \prod_{i,j} \gamma(s,\tau_i \times \pi_j,\psi)\]
And, if ${\bf G}_l$ is a classical group, then
\[ \gamma(s,\tau \times \pi,r,\psi)
   = \prod_{i=1}^e \gamma(s,\tau_i \times \pi_0,\psi) 
   \prod_{i=1}^e\prod_{j=1}^f 
   \gamma(s,\tau_i \times \pi_j,\psi) \gamma(s,\tau_i \times \tilde{\pi}_j,\psi). \]

Now, the second constituent of $R$ is given by
\begin{equation}\label{Siegelirred}
   r_2 = \left\{ \begin{array}{ll}
   		{\rm Sym}^2 \rho_m & \text{ if } {\bf G} = {\rm SO}_{2l+1} \\
		\wedge^2 \rho_m & \text{ if } {\bf G} = {\rm SO}_{2l} \text{ or } {\rm Sp}_{2l}.
	\end{array}. \right.
\end{equation}
In these cases, multiplicativity involves only $\tau$ of ${\rm GL}_m(F)$ given as a constituent of the induced representation of \eqref{gl:ind:rep}, and we have the following formula
\[ \gamma(s,\tau,r_2,\psi)
   = \prod_{i=1}^e \gamma(s,\tau_i,r_{m_i},\psi) 
   \prod_{i<j} \gamma(s,\tau_i \times \tau_j,\psi), \]
where $r_{m_i}$ is either ${\rm Sym}^2 \rho_{m_i}$ or $\wedge^2 \rho_{m_i}$, of the same kind as $r_2$.

\subsection{Global Langlands-Shahidi $L$-functions}

The global triples we consider $(\pi,r,\psi)$ require that $\pi = \otimes' \pi_v $ be globally generic. In particular, we can assume that $\pi$ is globally $\psi$-generic. By definition, there is a cusp form $\varphi$ in the space of $\pi$ such that
\begin{equation*}
   W_{M,\varphi}(m) = \int_{{\bf U}_M(k) \backslash {\bf U}_M(\mathbb{A}_k)} \varphi(um) \overline{\psi}(u) \, du \neq 0.
\end{equation*}

One extends $\varphi$ to an automorphic function
\[ \Phi : {\bf M}(k) {\bf U}(\mathbb{A}_k) \backslash {\bf G}(\mathbb{A}_k) \rightarrow \mathbb{C},
\]
as in \S~I.2.17 of \cite{MoWa1994}. We thus form a function in the space of the globally induced representation ${\rm I}(s,\pi)$, by setting
\begin{equation*}
   \Phi_s = \Phi \cdot \delta_{\!\alpha}^{s}, \ s \in \mathbb{C}.
\end{equation*}

The corresponding Eisenstein series is defined by
\begin{equation*}
   E(s,\Phi,g,{\bf P}) = \sum_{\gamma \in {\bf P}(k) \backslash {\bf G}(k)} \Phi_s(\gamma g),
\end{equation*}
which initially converges absolutely for $\Re(s) \gg 0$. For the globally generic representation $\pi$ of ${\bf M}(\mathbb{A}_k)$, the Fourier coefficient of the Eisenstein series $E(s,\Phi,g,{\bf P})$ is given by
\begin{equation*}
   E_\psi(s,\Phi,g,{\bf P}) = \int_{{\bf U}(k) \backslash {\bf U}(\mathbb{A}_k)} E(s,\Phi,ug,{\bf P}) \overline{\psi}(u) \, du.
\end{equation*}
Eisenstein series were studied in great generality by Langlands in \cite{La1976}, they have a meromorphic continuation and satisfy a functional equation. An argument of Harder \cite{Ha1974}, shows that in the case of function fields, Eisenstein series and their Fourier coefficients are rational functions on $q^{-s}$. The following can be found in \cite{ShFBook} for number fields, and in \cite{LoRationality} for function fields.

\begin{theorem}\label{rationalL}
Each Langlands-Shahid $L$-function $L(s,\pi,r_i)$ converges absolutely for $\Re(s) \gg 0$ and has a meromorphic continuation to the complex plane. If $k$ is a function field, then it is a rational function on $q^{-s}$.
\end{theorem}

In the previous theorem, we have completed $L$-functions $L(s,\pi,r_i)$. Going from partial $L$-functions $L^S(s,\pi,r_i)$, to completed ones, is a very delicate process. The Langlands-Shahidi method first makes the connection between Eisenstein series and Whittaker functionals
\begin{equation*}
   E_\psi(s,\Phi,g,{\bf P}) = \prod_v \lambda_{\psi_v}(s,\pi_v)({\rm I}(s,\pi_v)(g_v)f_{s,v}).
\end{equation*}
A careful computation at unramified places, then gives the connection to a product of partial $L$-functions
\begin{equation*}\label{eq1rationalL}
   E_\psi(s,\Phi,g,{\bf P}) = \prod_{v \in S} \lambda_{\psi_v}(s,\pi_v)({\rm I}(s,\pi_v)(g_v)f_{s,v}) \prod_{i=1}^{m(R)} L^S(1+is,\pi,r_i)^{-1}.
\end{equation*}
What we can obtain by arguing in this way, is a first form of the functional equation.

\begin{theorem}[Crude Functional Equation]\label{crudeFE}
Let $\pi$ be a globally $\psi$-generic cuspidal automorphic representation of ${\bf M}(\mathbb{A})$, let $R = \oplus_{i=1}^{m(R)} r_i$ be the adjoint representation. Then
   \begin{equation*}
      \prod_{i=1}^{m(R)} L^S(is,\pi,r_i) = \prod_{v \in S} C_{\psi}(s,\pi_v,\tilde{w}_0) \prod_{i=1}^{m(R)} L^S(1-is,\tilde{\pi},r_i).
   \end{equation*}
\end{theorem}

The induction result, Lemma~\ref{induction}, allows us to refine this into individual functional equations for each $i$, $1 \leq i \leq m(R)$, where we obtain local $\gamma$-factors at places $v$ in $S$:
\[ L^S(s,\pi,r_i) = \prod_{v \in S} \gamma(s,\pi_v,r_{i,v},\psi_v) 
   L^S(1-s,\tilde{\pi},r_i). \] 
Furthermore, from equation~\ref{temp:3}, and the discussion following it, we obtain $L$-functions and $\varepsilon$-factors at every non-archimedean place $v$, and for each $i$. At archimedean places, we refer to the crucial work of Shahidi \cite{ShF1985}. 

Since $\varepsilon$-factos are trivial for $v \notin S$, we can globally define
\[ \varepsilon(s,\pi,r_i) = \prod_v \varepsilon(s,\pi_v,r_{i,v},\psi_v). \]
We already know that partial $L$-functions initially converge on a right half plane, hence we now define completed $L$-functions
\[ L(s,\pi,r_i) = \prod_v L(s,\pi_v,r_{i,v}). \]
The Langlands-Shahidi method, thus produces the following form of the functional equation.

\begin{theorem}[Functional Equation]\label{LS:FE}
Let $(\pi,r,\psi)$ be a global Langlands-Shahidi triple, then
   \begin{equation*}
      L(s,\pi,r) = \varepsilon(s,\pi,r) L(1-s,\tilde{\pi},r).
   \end{equation*}
\end{theorem}

\section{Local Langlands}

\subsection{${\rm GL}_n$} Let $F$ be a non-archimedean local field and $n$ a positive integer. Let $\mathscr{A}_F(n)$ be the set of isomorphism classes of supercuspidal representations of ${\rm GL}_n(F)$ and let $\mathscr{G}_F(n)$ be the set of isomorphism classes of irreducible ${\rm Frob}$-semisimple $n$-dimensional representations of $\mathcal{W}_F'$.

\begin{theorem}[Local Langlands for ${\rm GL}_n$]\label{loc:L:GL}
For every non-archimedean local field $F$ and every integer $n \geq 1$, there is a correspondence
\begin{equation*}
   \mathcal{L}_{\rm loc} : \mathscr{A}_F(n) \rightarrow \mathscr{G}_F(n),
\end{equation*}
where we write $\sigma_\pi = \mathcal{L}_{\rm loc}(\pi)$. The maps are bijective, they satisfy and are uniquely characterized by the following properties:
\begin{enumerate}
   \item[(i)] $\mathcal{L}_{\rm loc}$ is the reciprocity map ${\rm Art}$ when $n=1$ and 
   \begin{equation*}
      \det(\sigma_\pi) = \omega_\pi \circ {\rm Art}.
   \end{equation*}
   \item[(ii)] If $\pi \in \mathscr{A}_F(n)$ and $\chi \in \mathscr{A}_F(1)$, we have that
   \begin{equation*}
      \mathcal{L}_{\rm loc}(\pi \cdot \chi) = \sigma_\pi \otimes \chi \text{ and }
   \tilde{\sigma}_\pi = \mathcal{L}_{\rm loc}(\tilde{\pi}). 
   \end{equation*}
   %\item[(i)] If $\pi$ is a subquotient of
%\begin{equation*}
   %{\rm Ind}(\pi_1 \otimes \cdots \otimes \pi_b),
%\end{equation*}
%with each $\pi_i$ irreducible quasi-tempered representation of ${\rm GL}_{n_i}(F)$, then
%\begin{equation*}
   %\sigma_\pi = \sigma_{\pi_1} \oplus \cdots \oplus \sigma_{\pi_b}.
%\end{equation*}
   \item[(iii)] $\mathcal{L}_{\rm loc}$ preserves and is uniquely determined by the property that it preserve local factors for every $\pi \in \mathscr{A}_F(m)$ and $\tau \in \mathscr{A}_F(n)$
\begin{align*}
   L(s,\pi \times \tau) &= L(s,\sigma_\pi \otimes \sigma_\tau) \\
   \varepsilon(s,\pi \times \tau,\psi) &= \varepsilon(s,\sigma_\pi \otimes \sigma_\tau,\psi),
\end{align*}
where $\psi: F \rightarrow \mathbb{C}^\times$ is a non-trivial character.
\end{enumerate}
\end{theorem}

The local Langlands conjecture was proved by Kutzko for ${\rm GL}_2$ \cite{Ku1980} followed by Henniart for ${\rm GL}_3$ \cite{He1983}. Over local function fields, it is a result of Laumon, Rapoport and Stuhler \cite{LaRaSt1993} for every $n$. For finite extensions of $\mathbb{Q}_p$, it was established independently by Harris-Taylor and Henniart \cite{HaTa2001,He2000}, and more recently by Scholze \cite{Sc2013}.

The uniqueness criterion presented in Theorem~\ref{loc:L:GL} using Rankin-Selberg local factors is due to Henniart \cite{He1993}, and originally appeared as an appendix to \cite{LaRaSt1993}; it is valid in any characteristic. However, we observe that Scholze's proof presents an alternative uniqueness criterion that makes the result of Harris-Taylor independent of Henniart's result. Namely, in \cite{Sc2013}, Scholze constructs functions $h$ and $f_{\tau,h}$ satisfying a trace relationship
\begin{equation*}
   {\rm Tr}(f_{\tau,h} \vert \pi) = {\rm Tr}(\tau \vert \mathcal{L}_{\rm loc}(\pi)) {\rm Tr}(h \vert \pi)
\end{equation*}
and proves that this condition uniquely characterizes the generalized reciprocity map $\mathcal{L}_{\rm loc}$. In addition, the map respects local $L$-functions and root numbers. The method of proof in this approach relies on powerful methods of Harris \cite{Ha1997} and Harris-Taylor \cite{HaTa2001}, based on Shimura varieties.

Local factors satisfy a multiplicativity property, where we observe that tempered representations for ${\rm GL}_n(F)$ are generic \cite{Ze1980}. Compatibility with the Langlands' classification, enables the passage from tempered to admissible representations.

\begin{enumerate}
   \item[(iv)] (Multiplicativity). \emph{Assume $\pi$ is the unique generic constituent
\begin{equation*}
   {\rm Ind}(\pi_1 \otimes \cdots \otimes \pi_b),
\end{equation*}
where each $\pi_i$ is a supercuspidal representation of ${\rm GL}_{n_i}(F)$, hence generic. Then we have the following multiplicativity property of local factors
\begin{align*}
   L(s,\pi \times \tau) &= \prod_{i,j} L(s,\pi_i \times \tau_j), \\
   \varepsilon(s,\pi \times \tau,\psi) &= \prod_{i,j} 
   \varepsilon(s,\pi_i \times \tau_j,\psi).
\end{align*}
}
   \item[(v)] (Langlands quotient). \emph{If $\pi$ is a subquotient of
\begin{equation*}
   {\rm Ind}(\pi_1 \otimes \cdots \otimes \pi_b),
\end{equation*}
with each $\pi_i$ irreducible quasi-tempered representation of ${\rm GL}_{n_i}(F)$, then
\begin{equation*}
   \sigma_\pi = \mathcal{L}_{\rm loc}(\pi) = 
   \sigma_{\pi_1} \oplus \cdots \oplus \sigma_{\pi_b}.
\end{equation*}
And we have compatibility with local $L$-functions and $\varepsilon$-factors for pairs.
}
\end{enumerate}

\subsection{Classical groups}
Given a classical group $\bf G$, if $F$ is a non-archimedean local field, we write $G = {\bf G}(F)$. Let $\mathscr{A}_F({\bf G})$ be the set of irreducible admissible representations of $G = {\bf G}(F)$. Let ${}^L\mathscr{G}_F({\bf G})$ be the set of $L$-parameters
\[ \phi : \mathcal{W}_F' \rightarrow {}^{L}G, \]
where the homomorphism $\phi$ is taken up to $\widehat{G}$ conjugacy for ${\rm SO}_{2n+1}$ and ${\rm Sp}_{2n}$, in our examples, cases (b) and (c), and up to ${\rm O}_{2n}(\mathbb{C})$-conjugacy in case (d).

An $L$-parameter $\phi \in {}^L\mathscr{G}_F({\bf G})$ is bounded if its image in ${}^LG$ projects onto a bounded subset of $\widehat{G}$. We say that $\phi \in {}^L\mathscr{G}_F({\bf G})$ is discrete if its image does not factor through a parabolic subgroup of ${}^LG$.

The Langlands correspondence for the classical groups in characteristic zero, both local and global, is due to Arthur \cite{Ar2013}, who provides an endoscopic classification via the trace formula, depending on the stabilization of the twisted trace formula proved by M\oe glin and Waldspurger in a series of papers culminating in \cite{MoWaX}. Under certain restrictions on the residual characteristic, the Langlands correspondence for local function fields was proved by Ganapathy-Varma \cite{GaVa2017}. Their method relies on the Kazhdan philosophy, on transferring information for representations of $\mathfrak{p}$-adic groups from characteristic $0$ to characteristic $p$, which is also established by Aubert, Baum, Plymen and Solleveld in \cite{AuBaPlSo2016}. We refer to \S~7 of \cite{GaLo2018}, for a description of the local Langlands correspondence over function fields for the classical groups and a general proof depending on a working hypothesis. The work of Genestier-Lafforgue \cite{GeLaPreprint}, now proves the local correspondence in characteristic $p$ for connected reductive groups.

\begin{theorem}\label{loc:L:classical}
For every split classical group $\bf G$ and every non-archimedean local field $F$ there is a finite-to-one map
\begin{equation*}
   \mathcal{L}_{\rm loc}: \mathscr{A}_F({\bf G}) 
   \rightarrow {}^L\mathscr{G}_F({\bf G}).
\end{equation*}
Writing $\phi_\pi = \mathcal{L}_{\rm loc}(\pi)$ for the corresponding Langlands parameter of a representation $\pi$, we have that
\begin{enumerate}
   \item[(i)] $\pi$ is a discrete series representation iff $\phi_\pi$ is a discrete $L$-parameter.
   \item[(ii)] $\pi$ is a tempered representation iff $\phi_\pi$ is a bounded $L$-parameter.
   \item[(iii)] Assume $\pi$ is the unique Langlands quotient of
   \begin{equation*}
      {\rm Ind}_P^G(\tau \cdot \eta),
   \end{equation*}
where $\tau$ is tempered and $\eta$ is in the Langlands situation. Then $\phi_\pi$ factors through ${}^LM$, i.e.,
\begin{equation*}
   \phi_\pi: \mathcal{W}_F' \xrightarrow{\phi_\tau} {}^LM_F \rightarrow {}^LG_F,
\end{equation*}
where $\phi_\tau = \mathcal{L}_{\rm loc}(\tau)$.
   \item[(iv)] Assume $\pi$ is a generic representation of $G$ and $\tau$ is a supercuspidal representation of ${\rm GL}_n(F)$, then
   \begin{align*}
      L(s,\pi \times \tau) &= L(s,\phi_\pi \otimes \phi_\tau) \\
      \varepsilon(s,\pi \times \tau,\psi) &= \varepsilon(s,\phi_\pi \otimes \phi_\tau,\psi).
   \end{align*}
   \item[(v)] For $\pi$ irreducible admissible non-generic and $\tau$ supercuspidal of ${\rm GL}_n(F)$, we have
   \begin{align*}
      \mu(s,\pi \times \tau,\psi) 
      &= \gamma(s,\phi_\pi^\vee \otimes \phi_\tau,\psi) \, 
      \gamma(-s,\phi_\pi \otimes \phi_\tau^\vee,\overline{\psi}) \\
      & \gamma(2s,r \circ \phi_\tau,\psi) \, \gamma(-2s,r^\vee \circ \phi_\tau,\overline{\psi}).
   \end{align*}
\end{enumerate}
Furthermore, $\mathcal{L}_{\rm loc}$ is characterized by properties $(i)$--$(v)$.
\end{theorem}

\subsection{On Local Langlands for reductive groups} The work of Arthur and M\oe glin-Waldpurger for the split classical groups was extended to quasi-split unitary groups by Mok \cite{Mo2015}, and recent work of Kaletha \cite{Ka2016} gives foundations for the general case of inner forms. Ongoing work of Fargues and Scholze aims towards the local Langlands conjecture when ${\rm char}(F) = 0$ and Genestier-Lafforgue have already established the supercuspidal case in ${\rm char}(F)=p$ \cite{GeLaPreprint}. The work of V. Lafforgue in characteristic $p$, and that of Scholze in characteristic zero, give the semisimple part of the $L$-parameter.

We only mention here that Base Change has played a crucial role in the local Langlands correspondence, starting with its proof for ${\rm GL}_n$. The article by Harris \cite{HaPreprint}, provides a Base Change approach towards local Langlands and related open questions.

\section{Global Langlands}

Let $k$ be a global field with Weil group $\mathcal{W}_k$. Let $\mathscr{A}_k(n)$ be the set of isomorphism classes of cuspidal automorphic representations of ${\rm GL}_n(\mathbb{A}_k)$ and let $\mathscr{G}_k(n)$ be the set of isomorphism classes of irreducible continuous degree $n$ representations of $\mathcal{W}_k$ which are unramified almost everywhere.

\begin{conjecture}[Global Langlands for ${\rm GL}_n$]\label{glob:L:GL}
For every $n \geq 2$, there is a correspondence
\begin{equation*}
   \mathcal{L}_{\rm glob}: \mathscr{A}_k(n) \rightarrow \mathscr{G}_k(n).
\end{equation*}
Writing $\Sigma = \mathcal{L}_{\rm glob}(\Pi)$, we have that
\begin{equation*}
   L(s,\Pi) = L(s,\Sigma)
\end{equation*}
and at every place $v$ of $k$ we have that
\begin{equation*}
   \Sigma_v = \mathcal{L}_{\rm loc}(\Pi_v).
\end{equation*}
The maps $\mathcal{L}_{\rm glob}$ are bijective and uniquely determined.
\end{conjecture}

The conjecture is known over function fields, where Drinfelds' Shtuka technique is the key to functoriality. L. Lafforgue makes a tour de force in his very important work for ${\rm GL}_n$ \cite{Dr1978}. However, it is an open problem over number fields. Cases where there are global representations $\Pi \in \mathscr{A}_k(n)$ and $\Sigma \in \mathscr{G}_k(n)$ which correspond to each other, are discussed in \S~\ref{globalization}.

\subsection{Converse Theorem} An important technique in proving that a representation of ${\rm GL}_n$ is automorphic is produced by the Converse Theorem of Cogdell and Piatetski-Shapiro \cite{CoPS1994}. For this we have a twisted version, where we fix a Gr\"o\ss encharakter character $\eta$ and a finite set $S$ of places of a global field $k$, including all archimedean places.

Given an integer $N$, consider the set $\mathcal{T}(S;\eta)$ consisting of $\tau = \tau_0 \otimes \eta$ such that: $\tau_0 \in \mathscr{A}_k(m)$, with $\tau_{0,v}$ is unramified for $ v \notin S$, and $m$ an integer ranging from $1 \leq m \leq N-1$.

\begin{theorem}\label{converse:thm}
Let $\Pi = \otimes' \Pi_v$ be an irreducible admissible representation of ${\rm GL}_N(\mathbb{A}_k)$ whose central character $\omega_\Pi$ is invariant under $k^{\times}$ and whose $L$-function $L(s,\Pi)$ is absolutely convergent for $\Re(s) \gg 0$. Suppose that for every $\tau \in \mathcal{T}(S;\eta)$ the $L$-function $L(s,\Pi \times \tau)$ is \emph{nice}, i.e., it satisfies:
\begin{enumerate}
\item[(i)] The functional equation
\begin{equation*}
   L(s,\Pi \times \tau) = \varepsilon(s,\Pi \times \tau)L(1-s,\widetilde{\Pi} \times \widetilde{\tau}).
\end{equation*}
\item[(ii)] $L(s,\Pi \times \tau)$ and $L(s,\widetilde{\Pi} \times \widetilde{\tau})$ are entire on $s \in \mathbb{C}$.
\item[(iii)] $L(s,\Pi \times \tau)$ and $L(s,\widetilde{\Pi} \times \widetilde{\tau})$ are bounded on vertical strips away from poles.
\end{enumerate} 
Then there exists an automorphic representation $\Pi'$ of ${\rm GL}_N(\mathbb{A}_k)$ such that $\Pi_v \cong \Pi'_v$ for all $v \notin S$.
\end{theorem}

\begin{remark}
If $k$ is a function field, we have the condition that $L(s,\Pi \times \tau)$ and $L(s,\widetilde{\Pi} \times \widetilde{\tau})$ be polynomials in the variable $q^{-s}$. This stronger condition, replaces conditions (ii) and (iii) in the above formulation of the Converse Theorem.
\end{remark}

\subsection{Functoriality for the classical groups} We restrict ourselves to generic representations \cite{CoKiPSSh2004,Lo2009}, since it leads to the Ramanujan Conjecture. We refer, to Arthur's trace formula approach to functoriality for the general case \cite{Ar2013}. A treatment of Arthur's trace formula and endoscopic classification would require a separate survey article.

We give an overview of the results of \cite{CoKiPSSh2004,Lo2009} on the globally generic functorial lift for the classical groups. Let ${\bf G}_n$ be a split classical group of rank $n$ defined over a global field $k$; $n \geq 2$ if ${\bf G}_n = {\rm SO}_{2n}$. Langlands functoriality is for ${\bf G}_n$ to ${\bf H}_N$, given by the following table.

\begin{center}
\begin{tabular}{|c|c|c|} \hline
   ${\bf G}_n$			& ${}^LG_n \hookrightarrow {}^LH_N$							& ${\bf H}_N$ \\ \hline
   ${\rm SO}_{2n+1}$	& ${\rm Sp}_{2n}(\mathbb{C}) \times \mathcal{W}_k \hookrightarrow {\rm GL}_{2n}(\mathbb{C}) \times \mathcal{W}_k$		& ${\rm GL}_{2n}$ \\
   ${\rm  Sp}_{2n}$		& ${\rm SO}_{2n+1}(\mathbb{C}) \times \mathcal{W}_k  \hookrightarrow {\rm GL}_{2n+1}(\mathbb{C}) \times \mathcal{W}_k$	& ${\rm GL}_{2n+1}$ \\
   ${\rm SO}_{2n}$		& ${\rm SO}_{2n}(\mathbb{C}) \times \mathcal{W}_k \hookrightarrow {\rm GL}_{2n}(\mathbb{C}) \times \mathcal{W}_k$		&${\rm GL}_{2n}$ \\
   \hline
\end{tabular}

\smallskip
Table~1
\end{center}

Locally, we let $\mathscr{A}_F^{\rm gen}({\bf G}_n)$ be the set of isomorphism classes of generic representations of ${\bf G}_n(F)$, when ${\bf G}_n$ is a classical group; and, we write $\mathscr{A}_F^{\rm gen}({\rm GL}_n)$, when we are in the case of a general linear group. Given Theorems~\ref{loc:L:GL} and \ref{loc:L:classical}, we have a map
\begin{equation*}
   \mathcal{L}_{\rm loc}^{\rm gen}: \mathscr{A}_F^{\rm gen}({\bf G}_n)
   \rightarrow \mathscr{A}_F^{\rm gen}({\rm GL}_N).
\end{equation*}
If we write $\Pi_0 = \mathcal{L}_{\rm loc}^{\rm gen}(\pi_0)$, the map is such that for every supercuspidal $\tau_0 \in \mathscr{A}_F^{\rm gen}({\rm GL}_m)$ and every non-trivial character $\psi_0: F \rightarrow \mathbb{C}^\times$ there is equality of $\gamma$-factors
\begin{equation*}
   \gamma(s,\Pi_0 \times \tau_0,\psi_0) = \gamma(s,\pi_0 \times \tau_0,\psi_0).
\end{equation*}
It is possible to show, with Henniart's uniqueness criterion \cite{He1993}, that preserving $\gamma$-factors as above determines the local Langlands lift. Furthermore, the map respects $L$-functions and $\varepsilon$-factors
\begin{align*}
   L(s,\Pi_0 \times \tau_0)& = L(s,\pi_0 \times \tau_0) \\
   \varepsilon(s,\Pi_0 \times \tau_0,\psi_0) &= \varepsilon(s,\pi_0 \times \tau_0,\psi_0).
\end{align*}

Globally, we let $\mathscr{A}_k^{\rm gen}({\bf G}_n)$ be the set of isomorphism classes of globally generic cuspidal automorphic representations of ${\bf G}_n(\mathbb{A}_k)$, when ${\bf G}_n$ is a classical group. We let $\mathscr{A}_k^{\rm gen}({\rm GL}_n)$ be the set of isomorphism classes of globally generic automorphic (not necessarily cuspidal) representations of ${\rm GL}_n(\mathbb{A}_k)$, when we are in the case of a general linear group.

\begin{theorem}[Classical groups]\label{fff}
There is a map
\begin{equation*}
   \mathcal{L}_{\rm gen}: \mathscr{A}_k^{\rm gen}({\bf G}_n) \rightarrow \mathscr{A}_k^{\rm gen}({\rm GL}_N).
\end{equation*}
Writing $\Pi = \mathcal{L}_{\rm gen}(\pi)$, we have that
\begin{equation*}
   L(s,\pi \times \tau) = L(s,\Pi \times \tau)
\end{equation*}
for every $\tau \in \mathscr{A}_k^{\rm gen}({\rm GL}_m)$. The representation $\Pi$ has trivial central character and can be expressed as an isobaric sum
\begin{equation*}
   \Pi = \Pi_1 \boxplus \cdots \boxplus \Pi_d,
\end{equation*}
where each $\Pi_i \in \mathscr{A}_k^{\rm gen}({\rm GL}_{N_i})$ is unitary self-dual cuspidal and $\Pi_i \ncong \Pi_j$ for $i \neq j$. Furthermore
\begin{equation*}
   \Pi_v = \mathcal{L}_{\rm loc}^{\rm gen}(\pi_v)
\end{equation*}
at every place $v$ of $k$.
\end{theorem}

Let us give a sketch of proof and then make a note on the image of functoriality. Let $\pi = \otimes' \pi_v \in \mathscr{A}_k^{\rm gen}({\bf G}_n)$. We construct a candidate lift $\Pi' = \otimes' \Pi'_v$, which will be a priori just admissible. We take $S$ to include all places where ramification may occur. %Then an application of the Converse Theorem will help us produce an automorphic representation $\Pi = \otimes' \Pi_v$, agreeing with $\Pi'$ at places outside of $S$. 

Initially, we understand the local Langlands correspondence for unramified representations $\pi_v$ of ${\bf G}_n(k_v)$ and $\Pi_v'$ of ${\rm GL}_N(k_v)$. Let these representations have Langlands parameters $\phi_{\pi_v}$ and $\phi_{\Pi_v'}$, they are then related by the following commutative diagram
\begin{center}
\begin{tikzpicture}
   \draw (1.95,0) node {$\mathcal{W}_{F}$};
   \draw[->,>=latex] (2.25,.25) -- (3.25,1);
   \draw (2.75,.4) node [right] {$\Phi_{\Pi_v'}$};
   \draw (0,1.25) node {${}^LG_n$};
   \draw[->,>=latex] (.5,1.25) -- (3.25,1.25);
   \draw (3.85,1.25) node {${}^L{\rm GL}_N$};
   \draw[->,>=latex] (1.5,.25) -- (.5,1);
   \draw (1,.4)  node[left] {$\phi_{\pi_v}$};
\end{tikzpicture}
\end{center}
possible via the Satake correspondence. Also, at archimedean places, we have compatibility with Artin parameters from the work of Shahidi \cite{ShFBook}.

The construction is such that global $L$-functions agree
\begin{equation*}
   L(s,\pi \times \tau) = L(s,\Pi' \times \tau)
\end{equation*}
for $\tau = \tau_0 \otimes \eta \in \mathcal{T}(S;\eta)$. And, at finite places $v \in S$, we have that $\tau = \tau_0 \otimes \eta \in \mathcal{T}(S;\eta)$ has $\tau_{0,v}$ unramified. Hence, $\tau_{0,v}$ for $v \in S$ is an irreducible constituent of an unramified principal series represnetation
\begin{equation*}
   {\rm Ind}(\chi_{1,v} \otimes \cdots \otimes \chi_{m,v}),
\end{equation*}
so that each $\chi_{i,v}: F^\times \rightarrow \mathbb{C}^\times$ is an unramified character. Then, we have
\begin{equation*}
   \gamma(s,\pi_v \otimes \tau_v,\psi_v) = 
   \prod_{i=1}^m \gamma(s,\pi_v \otimes (\chi_{i,v} \cdot \eta_v),\psi_v),
\end{equation*}
We can always take a Gr\"o\ss encharakter, such that $\eta_v$ is highly ramified at these places. Hence, the product on the right hand side of this last equation becomes stable.

By the Langlands-Shahidi method the global $L$-functions $L(s,\pi \times \tau)$ are nice, hence, so are $L(s,\Pi' \times \tau)$. An application of the Converse Theorem, then gives an automorphic representation $\Pi = \otimes \Pi_v$ of ${\rm GL}_N(\mathbb{A}_k)$ obtained from $\pi$ in such a way that $\Pi_v = \mathcal{L}_{\rm loc}^{\rm gen}(\pi_v)$ for $v \notin S$. From the classification of automorphic forms for ${\rm GL}_N$ \cite{JaSh1981}, $\Pi$ is the unique generic subquotient of a globally induced representation of the form ${\rm Ind}(\Pi_1 \otimes \cdots \otimes \Pi_d)$
with each $\Pi_i$ a cuspidal representation of ${\rm GL}_{N_i}(\mathbb{A}_k)$. We write this as an isobaric sum
\begin{equation*}
   \Pi =  \Pi_1 \boxplus \cdots \boxplus \Pi_d,
\end{equation*}
The representation $\Pi$ is unitary, what we need to show is that each $\Pi_i$ is unitary. For this, write $\Pi_i = \nu^{r_i} \Pi_{i,0}$ with each $\Pi_{i,0}$ unitary and $t_d \geq \cdots \geq t_1$, $t_1 \leq 0$. Because $\Pi$ is unitary, we cannot have $t_1>0$. Then
\begin{align*}
   L(s,\pi \times \widetilde{\Pi}_1) &= L(s,\Pi \times \widetilde{\Pi}_1) \\
   						&= \prod_{i} L(s + t_i,\Pi_i \times \widetilde{\Pi}_1). \nonumber
\end{align*}
An important analytic property, for Langlands-Shahidi $L$-functions for the classical groups, is that $L(s,\pi \times \tilde{\Pi}_1)$ is holomorphic for ${\rm Re}(s) > 1$. However, $L(s + t_1,\Pi_1 \times \tilde{\Pi}_1)$ has a simple pole at $s = 1 - t_1 \geq 1$. This pole carries through to a pole of $L(s,\pi \times \tilde{\Pi}_1)$, which gives a contradiction unless $t_1 = 0$. A recursive argument shows that all the $t_i$'s must be zero.

Also, every time we have $\Pi_i \cong \Pi_j$, we add to the multiplicity of the pole at $s=1$ of the $L$-function
\begin{equation}\label{eq1:fff}
   L(s,\pi \times \tilde{\Pi}_i) = \prod_{j} L(s,\Pi_j \times \tilde{\Pi}_i).
\end{equation}
However, $L(s,\pi \times \tilde{\Pi}_j)$ has at most a simple pole at $s = 1$. Hence, we must have $\Pi_i \ncong \Pi_j$, for $i \neq j$.

To show that each $\Pi_i$ is self-dual, we argue by contradiction and assume $\Pi_i \ncong \tilde{\Pi}_i$. The representation $\sigma = \tilde{\Pi}_i \otimes \tilde{\tau}$ of ${\bf M}_i(\mathbb{A}_k)$ satisfies $w_0(\sigma) \ncong \sigma$. Then, the $L$-function $L(s,\pi \times \tilde{\Pi}_i)$ has no poles. Now, we have that a pole appears on the right hand side of equation~\eqref{eq1:fff} from the $L$-function $L(s,\Pi_i \times \tilde{\Pi}_i)$. This is a contradiction, unless $\Pi_i$ is self-dual.

We now go through the classification of generic representations for the classical groups. First for supercuspidal representations, whic uses a globalization argument discussed in the following section. Then, there is the M\oe glin-Tadi\'c classification for discrete series. The standard module conjecture for tempered representations and the Langalnds classification. A crucial step is the stability property of $\gamma$-factors under twists by highly ramified characters.

\subsection{A note on the image of functoriality}

We observe that the argument found in \S~7 of \cite{LoRationality} is valid over any global field. That is, the argument provides a self conatined proof concerning the image of functoriality over number fleds of \cite{CoKiPSSh2004}, with the added note that symmetric and exterior square $L$-functions for cuspidal unitary automorphic representations are non-vanishing for $\Re(s)>1$, Proposition~\ref{nonvanishing}. For function fields, this non-vanishing is proved in \S~6 of \cite{LoRationality}.

\section{Globalization methods}\label{globalization}

Throughout this section assume that $F$ is a non-archimedean local field and $k$ is a global field such that $k_{v_0} \simeq F$ at a place $v_0$. We discuss globalization methods that take a representation $\pi$ over a non-archimedean local field $F$, for example, a supercuspidal representation of ${\rm GL}_n(F)$, and produce global representations $\Pi$ and $\Sigma$ that correspond to each other and $\Pi_{v_0} \simeq \pi$.

\begin{definition}
Let $\Pi$ be a cuspidal automorphic representation of ${\rm GL}_n(\mathbb{A}_k)$ and let $\Sigma$ be an $n$-dimensional $\ell$-adic representation of $\mathcal{W}_k$ on the Galois side. We say that the representations correspond to each other if
\begin{equation*}
   \Sigma = \mathcal{L}_{\rm glob}(\Pi),
\end{equation*}
as in Conjecture~\ref{glob:L:GL}.
\end{definition}

Note that if $\Pi$ and $\Sigma$ correspond to each other, then
\begin{equation*}
   \Sigma_v = \mathcal{L}_{\rm loc}(\Pi_v),
\end{equation*}
at every place $v$ of $k$.

\subsection{${\rm GL}_1$ and monomial representations for ${\rm GL}_n$} Already present in Jacquet-Langlands' work \cite{JaLa1970}, is the case of monomial representations, which provide a first family of cases where the global Langlands correspondence for ${\rm GL}_2$ is known. Locally, this arises from the tame case, and we have this for ${\rm GL}_n$, for example, when $p$ does not divide $n$. We begin with a basic lemma, which is a version of the Grundwald-Wang theorem of class field theory \cite{ArTa}.

\begin{lemma}\label{GL1:glob}
Let $\chi_0: {\rm GL}_1(F) \rightarrow \mathbb{C}^\times$ be a character. There exists a Gr\"o\ss encharakter $\chi : k^\times \backslash {\rm GL}_1(\mathbb{A}_k) \rightarrow \mathbb{C}^\times$ such that:
\begin{itemize}
   \item[(i)] $\chi_{v_0} \simeq \chi_0$;
   \item[(ii)] $\chi_v$ is unramified if $v \notin \left\{ v_0, \infty \right\}$;
   \item[(iii)] if ${\rm char}(F) = p$, then $\chi_\infty$ is tame.
\end{itemize}
\end{lemma}

Then, by global class field theory, there is a $1$-dimensional representation $\sigma = {\rm Art}(\chi)$ of $\mathcal{W}_k$ that corresponds to $\chi$.

From the theory of types \cite{BuKu1993}, we have level zero supercuspidal representations of ${\rm GL}_n(F)$. A level zero supercuspidal $\pi$ corresponds to a tamely ramified complex representation $\sigma$ of $\mathcal{W}_{F}$ via the local Langlands correspondence. More precisely, there exists an unramified extension $E/F$ of degree $n$ and a character $\eta_0:\mathcal{W}_E \rightarrow \mathbb{C}^\times$, such that
\begin{equation*}
   \sigma = {\rm Ind}_{\mathcal{W}_E}^{\mathcal{W}_F}(\eta_0).
\end{equation*}
For these representations we have the following globalization result.

\begin{proposition}
Let $\pi$ be a level zero supercuspidal representation of ${\rm GL}_n(F)$. Then there exist global representations $\Pi$ of ${\rm GL}_n(\mathbb{A}_k)$ and $\Sigma$ of $\mathcal{W}_k$ that correspond to each other and $\Pi_{v_0} \simeq \pi$.
\end{proposition}

To see this, we need the next basic lemma from Algebraic Number Theory.

\begin{lemma}
If $E/F$ is a finite separable extension, then there exists an extension $l/k$ and an $F$ isomorphism $E \simeq l \otimes_k k_v$.
\end{lemma}

Thus, given the extension $l/k$, an application Lemma~\ref{GL1:glob} gives a Gr\"o\ss encharakter $\eta : l^\times \backslash {\rm GL}_1(\mathbb{A}_l) \rightarrow \mathbb{C}^\times$ with $\eta_{v_0} \simeq \eta_0$. Then, we can set
\begin{equation*}
   \Sigma = {\rm Ind}_{\mathcal{W}_l}^{\mathcal{W}_k}(\eta).
\end{equation*}
Notice that we may choose, using Lemma~\ref{GL1:glob}, $\eta$ to be unramified at finite places distinct from $v_0$. Then at finite places $v$ in $S$ distinct from $v_0$, $\Pi_v$ is unramified. And, at finite unramified places we have the Satake parametrization \S~\ref{Lgroups:Satake}, hence
\begin{equation*}
   \Sigma_v = \mathcal{L}_{\rm loc}(\Pi_v),
\end{equation*}
at every place $v$ of $k$.

\subsection{On ${\rm GL}_n$ over number fields}\label{GL(n):0}
The global Langlands conjecture is an open problem over a number field $k$, even for ${\rm GL}_2$. However, there are important result of Langlands and Tunnell on base change, that allow us to obtain automorphic representations of ${\rm GL}_2(\mathbb{A}_k)$ that correspond to Galois representations.

We have the following globalizing lemma on the Galois side (Lemma~3.6 of \cite{He1983}, Lemma~4.13 of \cite{De1973}).

\begin{lemma}\label{GF:lemma}
Let $E/F$ be a finite Galois extension of non-archimedean local fields. Then there exist a global field $k$, a finite Galois extension $l/k$, with a place $v_0$ of $k$, and an isomorphism $\eta$ of $k_{v_0}$ onto of $F$ inducing an isomorphism of $k_{v_0} \otimes_k l$ onto $E$. In particular, $l$ has only one place $w_0$ above $v_0$ and the decomposition subgroup of ${\rm Gal}(l/k)$ at $w_0$ is itself, and identifies with ${\rm Gal}(E/F)$ via $\eta$.
\end{lemma}

A local Galois representation $\sigma$ of ${\rm Gal}(E/F)$, globalizes to a representation $\Sigma$ of ${\rm Gal}(l/k)$. Local Galois groups being solvable, the image of $\Sigma$ in the projective linear group is dihedral, $A_4$ or $S_4$ so there is by \cite{JaLa1970,La1980,Tu1978} a cuspidal automorphic representation $\Pi$ of ${\rm GL}_2(\mathbb{A}_l)$ so that
\[ \Sigma \longleftrightarrow \Pi, \]
correspond to each other.

For ${\rm GL}_n$, we have a globalization theorem of Harris \cite{Ha1997} and Harris-Taylor \cite{HaTa2001}, more recently refined by Scholze \cite{Sc2013} and Scholze-Shin \cite{ScSh2013}. These results, lead towards proofs of the local Langlands correspondence, by globalizing supercuspidal representations where there is a corresponding Galois representation. If we begin with a supercuspidal representation $\pi$ of ${\rm GL}_n(F)$, then there exists a number field $k$, with a place $v_0$ of $k$, and an automorphic representation $\Pi$ of ${\rm GL}_n(\mathbb{A})$ for which there is a global $n$-dimensional $\ell$-adic representation $\Sigma$ on the Galois side, such that
\[ \Pi_{v_0} \cong \pi \quad \text{and} \quad \Sigma_v = \mathcal{L}_{\rm loc}(\Pi_v), \]
at all places $v$ of $k$, i.e., the representations correspond to each other via the local Langlands correspondence for ${\rm GL}_n$. 

The work of many mathematicians aims to establish a large class of Galois representations which correspond to an automorphic representation under global Langlands. We refer to the work of Clozel and Thorne, starting with \cite{ClTh2014}, where several basic definitions can be found. And, the awe-inspiring 10 author effort \cite{ACC+} for potential automorphy.

Let $l$ be an imaginary CM field with a totally real subfield $k$, where $c \in {\rm Gal}(l/k)$ denotes conjugation. An automorphic representation $\Pi$ of ${\rm GL}_n(\mathbb{A}_l)$ is regular algebraic, conjugate self-dual, cuspidal (RACSDC), if for every place $v$ of $l$ such that $v \mid \infty$, the archimedean Langlands correspondence leads to a sum of pairwise distinct algebraic characters (regular algebraicity); $\Pi^c \cong \tilde{\Pi}$ (conjugate self-dual); and, $\Pi$ is cuspidal. Similar notions are RAECSDC, for essentially conjugate self-dual representations, and RAESDC, for essentially self-dual representations.

To summarize, such representations $\Pi$ of ${\rm GL}_n(\mathbb{A}_l)$ will have a corresponding global Galois representation after semisimplification
\[ \Pi \longleftrightarrow \Sigma, \]
such that
\begin{equation*}
   \Sigma_v = \mathcal{L}_{\rm loc}(\Pi_v),
\end{equation*}
at every place $v$ of $k$.

\subsection{Function fields}

The global field $k$ in characteristic $p$, arises as the field of functions of a smooth connected projective curve $X$ over the finite field $\mathbb{F}_q$. Each closed point of $X$, gives a valuation $v$ of $k$. With the \'etale topology, we obtain the the Galois group of $k$ from the fundamental group of $X$.

In the article \cite{Ka1986}, Katz shows that the morphism
\begin{equation*}
   {\rm Spec} \left(\mathbb{F}_q (\!( t^{-1} )\!) \right) \rightarrow 
   {\rm Spec}\left( \mathbb{F}_q \left[ t, t^{-1} \right] \right)
\end{equation*}
produces a functor
\begin{equation*}
   \left\{ \begin{array}{c} \text{finite `special' \'etale} \\
  				      \text{coverings of } 
  				      {\rm Spec}\left( \mathbb{F}_q \left[ t, t^{-1} \right] \right)
					\end{array} \right\}
   \longrightarrow 
   \left\{ \begin{array}{c} \text{finite \'etale coverings} \\
   				     \text{of } {\rm Spec} \left(\mathbb{F}_q (\!( t^{-1} )\!) \right)
					\end{array} \right\},
\end{equation*}
which is in fact an equivalence of categories. Here, a finite `special' \'etale cover is one of the punctured disk at two places $v_0$ and $v_\infty$, and is such that it is tame at $v_\infty$ and is unramifed for $v \notin \left\{ v_0, v_\infty \right\}$. On the right hand side, we have the fundamental group corresponding to the Galois group of $F$, which connects to the global Galois group of $k$ on the left hand side. If we begin with a local representation $\sigma_0$ on the right hand side, then the equivalence of categories proves the existence of a global Galois representation $\Sigma$, whose properties we summarize in Lemma~\ref{katz:loc-glob}.

We observe that when working with the \'etale topology, we need to pass between $\ell$-adic and complex representations. For this, we fix a prime number $\ell$, distinct from $p$, and an algebraic closure $\overline{\mathbb{Q}}_\ell$ of the field $\mathbb{Q}_\ell$ of $\ell$-adic numbers. We do this by fixing an isomorphism $\iota: \overline{\mathbb{Q}}_ \ell \simeq \mathbb{C}$.

\begin{lemma}\label{katz:loc-glob}
Let $\sigma_0$ be an irreducible $n$-dimensional Weil-Deligne representation of $\mathcal{W}_F$. Then, there is a rational function field $k$ with a set of two places $S = \left\{ v_0, v_\infty \right\}$, $k_{v_0} \simeq F$, and an irreducible continuous $n$-dimensional representation $\Sigma$ of $\mathcal{W}_k$ such that
\begin{itemize}
   \item[(i)] $\Sigma_{v_0} \simeq \sigma_0$;
   \item[(ii)] $\Sigma_v$ is unramified for $v \notin S$;
   \item[(iii)] $\Sigma_{v_\infty}$ is tame.
\end{itemize}
\end{lemma}

The result with prescribed ramification is for rational function fields. For a general function field, Katz proves that the representation can be globalized, yet there are finitely many other places where ramification can occur. Remark~\ref{Galois:loc-glob:rmk} below sketches a solution to the general problem with possibly one other tame place.

Now, working purely on the representation theoretic side, we obtain the corresponding globalization theorem for representations of ${\rm GL}_n$ in \cite{HeLo2013a}. The argument is generalized to a connected reductive group $\bf G$ in \cite{GaLo2018}, but introduces possibly finitely many tamely ramified places; it is refined in \cite{LoLS} to the following form.

\begin{lemma}\label{local-global:lemma}
Let $\pi$ be a supercuspidal representation of ${\bf G}(F)$. There exists a set of two places $S = \left\{ v_0, v_\infty \right\}$, $k_{v_0} \cong F$, and a cuspidal automorphic representation $\Pi = \otimes'_v \Pi_v$ of ${\bf G}(\mathbb{A}_k)$ satisfying the following properties:
\begin{itemize}
   \item[(i)] $\Pi_{v_0} \cong \pi$;
   \item[(ii)] $\Pi_v$ has an Iwahori fixed vector for $v \notin S$;
   \item[(iii)] $\Pi_{v_\infty}$ is a constituent of a tamely ramified principal series;
   \item[(iv)] if $\pi$ is generic, then $\Pi$ is globally generic. 
\end{itemize}
\end{lemma}

\begin{remark}\label{Galois:loc-glob:rmk}
The global Langlands correspondence for ${\bf G} = {\rm GL}_n$ over a function field $k$ \cite{LaL2002}, tells us that the cuspidal automorphic representation $\Pi$ of the previous Theorem corresponds to a Galois representation $\Sigma$. Because $\Pi$ and $\Sigma$ correspond to each other at every $v$, this extends Lemma~\ref{katz:loc-glob} on the Galois side to any $k$.
\end{remark}

\section{Ramanujan Conjecture}\label{Ramanujan}

For the group ${\rm GL}_n$, the character $\nu = \left| \det(\cdot) \right|$, obtained by composing the determinant with the absolute value, is defined for local and global fields. 

\subsection{${\rm GL}_n$}\label{RamanujanGL}
Let $F$ be a non-archimedean local field. An admissible unitary representation $\pi$ of ${\rm GL}_n(F)$ arises as a subquotient of a parabolically induced representation
\begin{equation}\label{localGLbound}
   \pi \hookrightarrow {\rm Ind}(\delta_0 \nu^{r_0} \otimes \cdots \otimes \delta_{b} \nu^{r_b}),
\end{equation}
where the $\delta_i$'s are discrete series representations of ${\rm GL}_{n_i}(F)$, $n = n_0 + \cdots + n_b$, and the Langlands parameters satisfy $0 \leq r_0 \leq \cdots \leq r_b$.

\smallskip

\noindent{\bf Ramanujan Conjecture for ${\rm GL}_n$.} \emph{Let $\Pi$ be a cuspidal automorphic representation of ${\rm GL}_n(\mathbb{A}_k)$, then each local component $\Pi_v$ is tempered. If $\Pi_v$ is unramified, temperedness implies that its Satake parameters satisfy
\begin{equation*}
   \left| \alpha_{j,v} \right|_{k_v} = 1, \quad 1 \leq j \leq n.
\end{equation*}
}

Over function fields, it is a theorem of Drinfeld for $n=2$ \cite{Dr1978} and L. Lafforgue for any $n$ \cite{LaL2002}. However, it is an open problem over number fields, even in the case of ${\rm GL}_2$ and $k = \mathbb{Q}$.

\subsection{Ramanujan bounds} Given a cuspidal (unitary) automorphic representation $\Pi = \otimes' \Pi_v$ of ${\rm GL}_n(\mathbb{A}_k)$, the Ramanujan Conjecture states that for every place $v$ of $k$, the representation $\Pi_v$ is isomorphic to a representation given by the form of equation \eqref{localGLbound} with $r_b = 0$; this exact bound is only known for a function field $k$.

A first bound that results from the study of $L$-functions for Rankin-Selberg products is what is known as the ``trival'' bound of $r_b \leq 1/2$ \cite{JaSh1981}. For arbitrary $n$ and any number field, it is Luo, Rudnick and Sarnak who established in \cite{LuRuSa1999} the bound
\begin{equation*}
   r_b \leq \dfrac{1}{2} - \dfrac{1}{n^2+1}.
\end{equation*}
It is interesting that this is still the best bound for general $n$, although breakthrough results have been obtained for $n \leq 4$.

There is the observation, made by Langlands that the automorphy of symmetric powers implies the Ramanujan Conjecture. The following conjecture is a very important open problem over number fields; it follows from the work of L. Lafforgue over function fields. 

\begin{conjecture}[Automorphy of symmetric powers]\label{sym:conj}
Given a cuspidal automorphic representation $\Pi$ of ${\rm GL}_n(\mathbb{A}_k)$, then the representation ${\rm Sym}^n(\Pi)$ is also automorphic for every $n$.
\end{conjecture}

The first breakthrough results in symmetric powers were obtained from the work of Kim and Shahidi, who use the Langlands-Shahidi method and establish the automorphy of symmetrice cube in \cite{KiSh2002}, in addition to the Ramanujan bound of $1/9$. Then, Kim carries their work further to establish symmetric fourth \cite{Ki2003} and with Sarnak \cite{KiSa2003}, they have the current best bound towards Ramanujan when $k = \mathbb{Q}$, namely
\begin{equation*}
   r_b \leq \dfrac{7}{64}.
\end{equation*}

Due to the presence or an infinite group of roots of unity in an arbitrary number field $k$, the matching result for ${\rm GL}_2(\mathbb{A}_k)$ had to wait until the work of Blomer and Brumley \cite{BlBr2011}. In addition to extending the Kim-Sarnak bound to arbitrary $k$, they obtain the record bounds of $r_b \leq 5/14$ for $n=3$ and $r_b \leq 9/22$ for $n=4$.

Blomer and Brumley, build on the result of Kim-Sarnak, using methods from Analytic Number Theory to study the series
\begin{equation*}
   L^S(s,\pi,{\rm Sym}^2) = \sum_{\mathfrak{a}} \lambda_{\pi,{\rm Sym}^2}(\mathfrak{a}) \mathcal{N}(\mathfrak{a})^{-s},
\end{equation*}
where the sum ranges over ideals $\mathfrak{a} \subset \mathcal{O}_k$ with $(\mathfrak{a},S)=1$, and the $\lambda_{\pi,{\rm Sym}^2}(\mathfrak{a})$'s are Fourier coefficients. They work with Deligne's bounds for hyper-Kloosterman sums to bound test functions and employ a Vonoroi summation formula with $S$-adic computations.

\begin{theorem}
If $L(s,\pi,{\rm Sym}^2)$ converges absolutely for $\Re(s) > 1$, then
\begin{equation*}
   r_b \leq \dfrac{1}{2} - \dfrac{1}{m+1}, \text{ for } m = \dfrac{n(n+1)}{2}.
\end{equation*}
\end{theorem}

We make a note, due to their appearence in several important places, concerning $L(s,\pi,{\rm Sym}^2)$ and the more general $L$-functions $L(s,\pi,{\rm Sym}^2 \boxtimes \chi)$.

\subsection{On twisted symmetric and exterior square $L$-functions}\label{symext:L}

In addition to their appearance towards the Ramanujan bounds, twisted symmetric and exterior square $L$-functions also play an important role in determining the image of functoriality. We make note that taking a trivial character suffices for the classical groups, non-trivial $\chi$ make an appearance for $\rm GSpin$ groups.

\begin{proposition}\label{nonvanishing}
Let $\tau$ be a cuspidal (unitary) automorphic representation of ${\rm GL}_n(\mathbb{A}_k)$ and $\chi$ a unitary Gr\"o\ss encharakter. Let $r = {\rm Sym}^2$ or $\wedge^2$. Then $L(s,\tau,r \boxtimes \chi)$ is holomorphic except for possible poles appearing at $s=0$ or $s=1$; it is non-vanishing for $\Re(s)>1$.
\end{proposition}

Twisted symmetric and exterior square $L$-functions are available both via an integral representation and the Langlands-Shahidi method. However, the proof of this proposition was recently completed by Takeda \cite{TaS2014,TaS2015}, who builds upon results of Bump and Ginzburg \cite{BuGi1992}, using an integral representation for the partial $L$-functions $L^S(s,\tau,{\rm Sym}^2\boxtimes\chi)$. It is instructive to recall that non-vanishing is obtained first for the Rankin-Selberg $L$-function. Indeed, Corollary~5.8 of \cite{JaSh1981} states that $L(s,\tau \times (\tau \cdot \chi))$, which is holomorphic and non-vanishing for $\Re(s)>1$. We then use the relationship
\begin{equation*}
   L(s,\tau \times (\tau \cdot \chi)) = 
   L(s,\tau,{\rm Sym}^2) L(s,\tau,\wedge^2).
\end{equation*}
This well known relationship is proved by a local-to-global argument, where we use \cite{ShFBook} over number fields and Lemma~\ref{local-global:lemma} over function fields.

We observe that exterior square $L$-functions were studied by Jacquet-Shalika, who produce an integral representation in \cite{JaSh1990}. Furthermore, it is possible to show that the integral representation approach and that of the Langlands-Shahidi method are compatible. While Takeda's result is for partial $L$-functions, at the remaining possibly ramified places, local $L$-functions are non-zero and holomorphic for $\Re(s)>1$. To see this, we use the tempered $L$-function conjecture, proved in \cite{HeOp2013}, which gives holomorphy of local $L$-functions for $\Re(s)>0$ if $\pi$ is tempered, then, Langlands classification and the trivial Ramanujan bound $r_b \leq 1/2$ for ${\rm GL}_n$. %and, for the classical groups, we will get this for $\Re(s)>1$ using the unitary generic bound of Lapid-Mui\'c-Tadi\'c--).

\subsection{Cases towards Ramanujan for ${\rm GL}_n$}

The globalization techinques discussed in \S~\ref{GL(n):0}, give cases where the Ramanujan Conjecture holds to be true. There is work of Shin \cite{Shin2011} and Caraiani \cite{Ca2012}; the series of papers by Clozel-Thorne, starting with \cite{ClTh2014}, which have applications up to symmetric $8^{\rm th}$ powers; and, the current 10 author work in progress \cite{ACC+}.

\begin{theorem}
If $\Pi$ is a cuspidal automorphic form of ${\rm GL}_n(\mathbb{A}_k)$ constructed as in \S~\ref{GL(n):0} and has corresponding Galois representation $\Sigma$ of $\mathcal{W}_k$. Then $\Pi$ satisfies the Ramanujan Conjecture.
\end{theorem}

\subsection{{\bf G} a reductive group}\label{GenRamanujan}
If $\bf G$ is a connected reductive group over a global field $k$, then the Ramanujan conjecture no longer holds for all cuspidal automorphic representations $\pi$ of ${\bf G}(\mathbb{A}_k)$. Counterexamples were produced by Howe and Piatetski-Shapiro for ${\bf G} = {\rm Sp}_4$ and ${\rm U}_3$, via the $\theta$-correspondence \cite{HoPS1979}. However, if one adds the genericity condition, then it is expected to be true.%, see for example Theorem~\ref{classical:Ram}.

\subsection*{Generalized Ramanujan Conjecture}\label{GRamanujan}
\emph{Let $\Pi = \otimes' \Pi_v$ be a globally generic cuspidal automorphic representation of ${\bf G}(\mathbb{A}_k)$. Then each $\Pi_v$ is tempered. If $\Pi_v$ is unramified, then its Satake parameters satisfy
\begin{equation*}
   \left| \alpha_{j,v} \right|_{k_v} = 1, \quad 1 \leq j \leq n.
\end{equation*}
}

\smallskip

The globally generic lift from the classical groups to ${\rm GL}_N$ provides us with a way to establish new cases of the Generalized Ramanujan Conjecture over function fields thanks to the exact Ramanujan bounds for ${\rm GL}_N$ of L. Lafforgue. Over a number field, we can import the best known bounds for ${\rm GL}_N$, discussed in the previous section, to the classical groups \cite{Ar2013,CoKiPSSh2004}. Over function fields, we have the Ramanujan conjecture for the classical and unitary groups over funtion fields \cite{Lo2009,LoLS}.

\begin{theorem}\label{classical:Ram}
Let $\bf G$ be either ${\rm SO}_{2n+1}$, ${\rm Sp}_{2n}$ or ${\rm SO}_{2n}$ defined over a function field $k$. Then, a globally generic cuspidal automorphic representation $\Pi$ of ${\bf G}_n(\mathbb{A}_k)$, satisfies the Ramanujan Conjecture.
\end{theorem}

For a general connected reductive group over a function field $k$, V. Lafforgue has established a very important automorphic to Galois lift, associating a representation $\Sigma$ to a cuspidal automorphic representation $\Pi$ of ${\bf G}(\mathbb{A}_k)$. This construction agrees at every unramified place with the Satake parametrization. We also note recent work of Sawin and Templier, who prove under certain conditions a general Ramanujan conjecture over function fields \cite{SaTePreprint}.

%Combining the Langlands-Shahidi method over function fields with the results of L. Lafforgue and V. Lafforgue, makes possible the following result; which is non vacuous, due to known cases of the Ramanujan Conjecture over function fields \cite{LoLS}.

%\begin{theorem}\label{RamRie:thm}
%Let $(\pi,r,\psi)$ be a global Langlands-Shahidi triple such that $r$ is faithful, assume that $\pi$ satisfies the Ramanujan Conjecture. Then, there is an automorphic representation of ${\rm GL}_N(\mathbb{A}_k)$
%\begin{equation}\label{eq:isobaric}
   %\mathcal{T} = \mathcal{T}_1 \boxplus \cdots \boxplus \mathcal{T}_d,
%\end{equation}
%written as an isobaric sum of unitary cuspidal automorphic representations $\mathcal{T}_j$ of ${\rm GL}_{N_j}(\mathbb{A}_k)$, such that
%\begin{equation*}
   %L^S(s,\pi,r) = L^S(s,\mathcal{T})
%\end{equation*}
%and at every place $v$ of $k$ we have that
%\begin{equation*}
   %\gamma(s,\pi_v,r_{v},\psi_v) = \gamma(s,\mathcal{T}_v,\psi_v).
%\end{equation*}
%Furthermore, the zeros of $L(s,\pi,r)$ are contained in the line $\Re(s) = 1/2$.
%\end{theorem}


\begin{thebibliography}{99}

\bibitem{ACC+} P. Allen, F. Calegari, A. Caraiani, T. Gee, D. Helm, B. Le Hung, H. Newton, P. Scholze, R. Taylor, and J. Thorne, \emph{Potential automorphy over CM fields}, work in progress.

\bibitem{Ar2003} J. Arthur, \emph{The principle of functoriality}, Bull. A.M.S. {\bf 40} (2003), 39-53.

\bibitem{Ar2013} \rule{3em}{.5pt}\thinspace, \emph{The Endoscopic Classification of Representations: Orthogonal and Symplectic Groups}, Colloquium Publications of the AMS, Vol. 61, 2013, Providence, RI.

\bibitem{Au2018} A.-M. Aubert, \emph{Around the Langlands Program}, Jahresber Dtsch Math-Ver {\bf 120} (2018), 3-40.

\bibitem{ArTa} E. Artin and J. Tate, \emph{Class field theory}, Amer. Math. Soc., Providence, RI, 2009.

\bibitem{AuBaPlSo2016} A.-M. Aubert, P. Baum, R. Plymen and M. Solleveld,  \emph{The local Langlands correspondence for inner forms of ${\rm SL}_n$}, Res. Math. Sci. (2016) 3-32.

\bibitem{BeGe2004} J. Bernstein and S. Gelbart (Eds.), \emph{An introduction to the Langlands program}, Birkh\"auser, 2004.

\bibitem{BlBr2011} V. Blomer and F. Brumley, \emph{On the Ramanujan conjecture over number fields}, Ann. Math {\bf } (2011), 581-605. 

\bibitem{BlBr2013} \rule{3em}{.5pt}\thinspace, \emph{The role of the Ramanujan Conjecture in analytic number theory}, Bull. A.M.S. {\bf 50} (2013), 267-320.

\bibitem{Bo1979} A. Borel, \emph{Automorphic $L$-functions}, Proc. Symp. Pure Math., vol. 33, part 2, Amer. Math. Soc., Providence, RI, 1979, pp. 27-61. 

\bibitem{BoCa1979} A. Borel and W. Casselman (editors), \emph{Automorphic Forms, Representations, and $L$-functions}, Proc. Symp. Pure Math., vol. 33, parts 1 \& 2, Amer. Math. Soc., Providence, RI, 1979

\bibitem{BoWa} A. Borel and N. Wallach, \emph{Continuous cohomology, discrete subgroups, and representations of reductive groups},
2nd ed., Math. Surv. and Mono. {\bf 67}, AMS, 2000.

\bibitem{Bou} N. Bourbaki, \emph{Groupes et alg\`ebres de Lie}, chapitres 4, 5 et 6, Hermann, Paris 1968.

\bibitem{BuGi1992} D. Bump and D. Ginzburg, \emph{Symmetric square $L$-functions on ${\rm GL}(r)$}, Ann. Math. {\bf 136} (1992), 137-205.

\bibitem{BuHe2006} C. J.  Bushnell and G. Henniart, \emph{The local Langlands conjecture for ${\rm GL}(2)$}, Grund. Math. Wiss. {\bf 335}, Springer-Verlang, Berlin, 2006.

\bibitem{BuKu1993} C. J. Bushnell and P. C. Kutzko, \emph{The admissible dual of ${\rm GL}(N)$ via compact open subgroups}, Ann. Math. Studies, {\bf 129}, Princeton University Press, 1998.

\bibitem{Ca2012} A. Caraiani, \emph{Local-global compatibility and the action of monodromy on nearby cycles}, Duke Math. J. {\bf 161} (2012), 2311-2413.

\bibitem{Ca1979} P. Cartier, \emph{Representations of $\mathfrak{p}$-adic groups: a survey}, Proc. Symp. Pure Math. {\bf 33} (1979), vol. 1, 111-156.

\bibitem{CaNotes} W. Casselman, \emph{Introduction to the theory of admissible representations of $p$-adic reductive groups}, unpublished notes.

\bibitem{CaFr} J. W. S. Cassels and A. Fr\"olich, editors, \emph{Algebraic Number Theory}, $2^{\rm nd}$ ed., London Math. Soc., London, 2010.

\bibitem{ClTh2014} L. Clozel and J. A. Thorne, \emph{Level-raising and symmetric power functoriality, I}, Comp. Math. {\bf 150} (2014), 729-748.

\bibitem{CoKiPSSh2004} \rule{3em}{.5pt}\thinspace, \emph{Functoriality for the classical groups}, Publ. Math. IHES {\bf 99} (2004), 163-233.

\bibitem{CoPS1994} J. W. Cogdell and I. I. Piatetski-Shapiro, \emph{Converse theorems for ${\rm GL}_n$}, Publ. Math. IH\'ES {\bf 79} (1994), 157-214.

\bibitem{CoPS2004} \rule{3em}{.5pt}\thinspace, \emph{Remarks on Rankin-Selberg convolutions}, (published in [H. Hida, D. Ramakrishnan and F. Shahidi, editors, Contributions to Automorphic Forms, Geometry and Number Theory, John Hopkins University Press, 2004, 255-278]).

\bibitem{De1973} P. Deligne, \emph{Les constantes des \'equations fonctionnelles des fonctions $L$}, International Summer School on Modular Functions, Antwerp II, Lecture Notes in Math., vol. 349, Springer-Verlag, 1973, pp. 501-597.

\bibitem{De1974} \rule{4.3em}{.5pt}\thinspace, \emph{La conjecture de Weil. I \& II}, Pub. Math. IH\'ES {\bf 43} (1974), 273-307, \& {\bf 52} (1980), 137-252.

\bibitem{SGA3} M. Demazure and A. Grothendieck, \emph{Sch\'emas en Groupes I, II \& III}, S\'eminaire de G\'eometrie Alg\'ebrique du Bois Marie 1962/64, SGA 3, Lecture Notes in Math., vols. 151, 152 \& 153, Springer-Verlag, 1970.

\bibitem{Dr1978} V. G. Drinfeld, \emph{Langlands' conjecture for ${\rm GL}(2)$ over function fields}, Proceedings of the International Congress of Mathematicians, Helsinki, 1978.

\bibitem{Dr1988} \rule{4.3em}{.5pt}\thinspace, \emph{The proof of Petersson's conjecture for ${\rm GL}(2)$ over a global field of characteristic $p$}, Functional Analysis and its Applications {\bf 22} (1988), 23-43.

\bibitem{Fl1979} D. Flath, \emph{Decomposition of representations into tensor products}, Proc. Symp. Pure Math., vol. 33, part 2, Amer. Math. Soc., Providence, RI, 1979, pp. 179-184.

\bibitem{Fr2013} E. Frenkel, \emph{Langlands Program, Trace Formulas, and their Geometrization}, Bull. A.M.S. {\bf 50}, 1-55.

\bibitem{GaLo2018} W. T. Gan and L. Lomel\'i, \emph{Globalization of supercuspidals over function fields and applications}, J. European Math. Soc. {\bf 20} (2018), 2813-2858.

\bibitem{Ga2015} R. Ganapathy, \emph{The local Langlands correspondence for ${ \rm GSp}_4$ over local function fields}, Amer. J. Math. {\bf 137} (2015), 1441-1534.

\bibitem{GaVa2017} R. Ganapathy and S. Varma, \emph{On the local Langlands correspondence for the split classical groups over local function fields}, J. Inst. Math. Jussieu {\bf 16} (2017), 987-1074.

\bibitem{Ge1984} S. Gelbart, \emph{An elementary introduction to the Langlands Program}, Bull. A.M.S. {\bf 10}, 177-219.

\bibitem{GeLaPreprint} A. Genestier and V. Lafforgue, \emph{Chtoucas restreints pour les groupes r\'eductifs et param\'etrisation de Langlands locale}, preprint.

\bibitem{GoJa1972} R. Godement and H. Jacquet, \emph{Zeta functions of simple algebras}, Lecture Notes in Math., vol. 260, Springer-Verlag, 1972.

\bibitem{HaSatake} T. J. Haines, \emph{On Satake parameters for representations with parahoric fixed vectors}, preprint.

\bibitem{Ha1974} G. Harder, \emph{Chevalley groups over function fields and automorphic forms}, Annals of Math. {\bf 100} (1974), 249-306.

\bibitem{Ha1997} M. Harris, \emph{Supercuspidal representations in the cohomology of Drinfel'd upper half spaces; elaboration of Carayol's program}, Invent. Math. {\bf 129} (1997), 75-119.

\bibitem{HaPreprint} \rule{4.3em}{.5pt}\thinspace, \emph{Incorrigible representations}, preprint.

\texttt{http://arxiv.org/abs/1811.05050}.

\bibitem{HaTa2001} M. Harris and R. Taylor, \emph{The geometry and cohomology of some simple Shimura varieties}, Annals of Mathematics Studies, vol. 151, Princeton University Press, Princeton, NJ, 2001, With an appendix by Vladimir G. Berkovich.

\bibitem{HeOp2013} V. Heiermann and E. Opdam, \emph{On the tempered $L$-function conjecture}, Amer. J. Math. {\bf 137} (2013), 777-800.

\bibitem{He1983} G. Henniart, \emph{La conjecture de Langlands locale pour ${\rm GL}(3)$}, M\'em. Soc. Math. France {\bf 11 -12} (1983), 1-186.

\bibitem{He1993} \rule{3em}{.5pt}\thinspace, \emph{Caract\'erisation de la correspondance de Langlands locale par les facteurs $\varepsilon$ de paires}, Invent. Math. {\bf 113} (1993), 339-350.

\bibitem{He2000} \rule{3em}{.5pt}\thinspace, \emph{Une preuve simple des conjectures de Langlands pour ${\rm GL}(n)$ sur un corps $p$-adique}, Invent. Math. \textbf{139} (2000), no.~2, 439--456.

\bibitem{HeLo2013a} G. Henniart and L. Lomel\'i, \emph{Uniqueness of Rankin-Selberg factors}, Journal of Number Theory {\bf 133} (2013), 4024-4035.

\bibitem{HoPS1979} R. Howe and I. Piatetski-Shapiro, \emph{A counterexample to the `generalized Ramanujan conjecture' for (quasi)-split groups'}. (Published in [Automorphic Forms, Representations and $L$-functions, Proc. Symp. Pure Math. {\bf 33}, part I, AMS, Providence, 1979, pp. 315-322]).

\bibitem{JaLa1970} H. Jacquet and R. P. Langlands, \emph{Automorphic Forms on ${\rm GL}(2)$}, Lecture Notes in Math. {\bf 114}, Springer-Verlag, Berlin 1970.

\bibitem{JaPSSh1983} H. Jacquet, I. I. Piatetski-Shapiro and J. A. Shalika, \emph{Rankin-Selberg convolutions}, Amer. J. Math. {\bf 105} (1983), 387-464.

\bibitem{JaSh1981} H. Jacquet and J. A. Shalika, \emph{On Euler products and the classification of automorphic representations I, II},
Am. J. Math. {\bf 103} (1981), 499-558 and 777-815.

\bibitem{JaSh1990} \rule{3em}{.5pt}\thinspace, \emph{Exterior square $L$-functions}, in [Automorphic forms, Shimura varieties, and $L$-functions, vol. II (Ann Arbor, MI, 1988), Perspect. Math. {\bf 11}, Academic Press, Boston, MA, 1990, pp. 143-226].

\bibitem{Ka2016} T. Kaletha, \emph{Rigid inner forms of real and $p$-adic groups}, Ann. Math. {\bf 2016}, 559-632.

\bibitem{Ka1986} N. M. Katz, \emph{Local-to-global extensions of representations of fundamental groups}, Ann. Inst. Fourier ${\bf 1986}$, 69-106.

\bibitem{KeSh1988} C. D. Keys and F. Shahidi, \emph{Artin $L$-functions and normalization of intertwining operatons}, Ann. Scient. \'E.N.S. {\bf 21} (1998), 67-89.

\bibitem{Ki2003} H. H. Kim, \emph{Functoriality for the exterior square of ${\rm GL}_4$ and the symmetric fourth of ${\rm GL}_2$}, J. A. M. S. {\bf 16} (2003), 139-183.

\bibitem{KiSa2003} H. H. Kim and P. Sarnak, \emph{Refined estimates towards the Ramanujan and Selberg Conjectures}, J. Amer. Math. Soc. {\bf 16} (2003), 139-183, Apendix to H. Kim, \emph{Functoriality for the exterior square of ${\rm GL}(4)$ and symmetric fourth of ${\rm GL}(2)$}.

\bibitem{KiSh2002} H. H. Kim and F. Shahidi, \emph{Functorial products for ${\rm GL}_2 \times {\rm GL}_3$ and the symmetric cube for ${\rm GL}_2$}, Ann. Math. {\bf 155} (2002), 837-893.

\bibitem{KiJL2007} J. L. Kim, \emph{Supercuspidal representations: an exhaustion theorem}, J.A.M.S. {\bf 20} (2007), 273-320.

\bibitem{Kn1991} M.-A. Knus, \emph{Quadratic and hermitian forms over rings}, Springer, Berlin, 1991.

\bibitem{Ku1980} P. Kutzko, \emph{The Langlands conjecture for ${\rm GL}(2)$ of a local field}, Ann. Math. {\bf 112} (1980), 381-412.

\bibitem{LaL2002} L. Lafforgue, \emph{Chtoucas de Drinfeld et correspondance de Langlands}, Invent. Math. {\bf 147} (2002), 1-241.

\bibitem{LaV} V. Lafforgue, \emph{Chtoucas pour les groupes r\'eductifs et param\'etrization de Langlands globale}, J.A.M.S. {\bf 31} (2018), 719-891.

\bibitem{La1967} R. P. Langlands, \emph{Euler Products}, Yale University Press, New Haven, 1967.

\bibitem{La1976} \rule{3em}{.5pt}\thinspace, \emph{On the functional equations satisfied by Eisenstein series}, Lecture Notes Math. {\bf 544}, Springer, Berlin, 1976.

\bibitem{LaCorvallisI} \rule{3em}{.5pt}\thinspace, \emph{On the notion of an automorphic representation}, Proc. Symp. Pure Math., {\bf 33}, part 1 (1979), 203-207.

\bibitem{LaCorvallisII} \rule{3em}{.5pt}\thinspace, \emph{Automorphic representations, Shimura varieties, and motives. Ein M\"archen}, Proc. Symp. Pure Math., {\bf 33}, part 2 (1979), 205-246.

\bibitem{LaYale} \rule{3em}{.5pt}\thinspace, \emph{On the functional equation of Artin's $L$-functions}, mimeographed notes, Yale University.

\bibitem{La1980} \rule{3em}{.5pt}\thinspace, \emph{Base Change for ${\rm GL}(2)$}, Ann. Math. Studies, vol. 96, Princeton Univ. Press, New Jersey, 1980.

\bibitem{LaMuTa2004} E. Lapid, G. Mui\'c and M. Tadi\'c, \emph{On the generic unitary dual of quasisplit classical groups}, Internat. Math. Res. Not. {\bf 26} (2004), 1335-1354.

\bibitem{LaRaSt1993} G. Laumon, M. Rapoport, and U. Stuhler, \emph{$\mathscr{D}$-elliptic sheaves and the Langlands correspondence}, Invent. Math. {\bf 113} (1993) 217-338.

\bibitem{Lo2009} L. Lomel\'i, \emph{Functoriality for the classical groups over function fields}, Internat. Math. Res. Notices {\bf 22} (2009), 4271-4335.

\bibitem{LoRationality} \rule{3em}{.5pt}\thinspace, \emph{Rationality and holomorphy of Langlands-Shahidi $L$-functions over function fields}, to appear in Math. Zeit. (2018). https://doi.org/10.1007/s00209-018-2100-7

\bibitem{LoLS} \rule{3em}{.5pt}\thinspace, \emph{The Langlands-Shahidi method over function fields},

\texttt{http://arxiv.org/abs/1507.03625}.

%\bibitem{LoRam} \rule{3em}{.5pt}\thinspace, \emph{Ramanujan conjecture and Riemann Hypothesis for the unitary groups},

%\texttt{http://arxiv.org/abs/1507.03625}.

\bibitem{LuRuSa1999} W. Z. Luo, Z. Rudnick and P. C. Sarnak, \emph{On the generalized Ramanujan conjecture for ${\rm GL}(n)$}, Proc. Symp. Pure Math. {\bf 66}, part 2 (1999), 301-310.

\bibitem{Ma1958} F. I. Mautner, \emph{Spherical functions over $\mathfrak{P}$-adic fields, I}, Amer. J. Math {\bf 80} (1958), 441-457.

\bibitem{Ma1964} \rule{3em}{.5pt}\thinspace, \emph{Spherical functions over $\mathfrak{P}$-adic fields, II}, Amer. J. Math {\bf 86} (1964), 171-200.

\bibitem{MiBook} J. Milne, \emph{Algebraic group schemes}, 2017.

\bibitem{MoWa1994} C. M\oe glin and J.-L. Waldspurger, \emph{D\'ecomposition spectrale et s\'eries d'Eisenstein: une paraphrase de l'\'ecriture}, Progress in Math. {\bf 113}, Birkh\"auser Verlag, Basel 1994.

\bibitem{MoWaX} \rule{3em}{.5pt}\thinspace, \emph{Stabilisation de la formule des traces tordue X: stabilisation spectrale}, preprint.

\texttt{http://arxiv.org/abs/1412.2981}

\bibitem{Mo2015} C. P. Mok, \emph{Endoscopic classification of representations of quasi-split unitary groups}, Memoirs AMS {\bf 235} (2015), 1-248. 

\bibitem{Mo1982} L. Morris, \emph{Eisenstein series for reductive groups over global function fields I, II},
Can. J. Math. {\bf 34} (1982), 91-168 and 1112-1182.

\bibitem{RaVa1999} D. Ramakrishnan and R. J. Valenza, \emph{Fourier Analysis on Number Fields}, Graduate Texts in Math. {\bf 186}, Springer, New York, 1999.

\bibitem{Sa2005} P. Sarnak, \emph{Notes on the Generalized Ramanujan Conjectures}, Clay Math. Proc. {\bf 4} (2005), 659-685.

\bibitem{Sa1963} I. Satake, \emph{Theory of spherical functions on reductive algebraic groups over $p$-adic fields}, Publ. Math. I.H.\'E.S {\bf 18} (1963), 5-69.

\bibitem{SaTePreprint} W. Sawin and N. Templier, \emph{On the Ramanujan conjecture for automorphic forms over function fields I. Geometry}, preprint.

\texttt{http://arxiv.org/abs/1805.12231}

\bibitem{Sc2013} P. Scholze, \emph{The local Langlands correspondence for ${\rm GL}_n$ over $p$-adic fields}, Invent. Math. {\bf 192} (2013), 663-715.

\bibitem{ScSh2013} P. Scholze and S. W. Shin, \emph{On the cohomology of compact unitary group shimura varieties at ramified split places}, J.A.M.S. {\bf 26} (2013), 261-294.

\bibitem{Se} J.-P. Serre, {Cours d'arithm\'etique}, $4$\`eme \'ed., Puf, Paris, 1994.

\bibitem{ShFBook} F. Shahidi, \emph{Eisenstein series and Automorphic $L$-functions}, Collo. Pub., vol. 58, AMS, Providence, RI, 2010.

\bibitem{ShF2004} \rule{3em}{.5pt}\thinspace, \emph{On the Ramanujan Conjecture for quasisplit groups}, Asian J. Math {\bf 8} (2004), 813-836.

\bibitem{ShF1985} \rule{3em}{.5pt}\thinspace, \emph{Local coefficients as Artin factors for real groups}, Duke Math. J. {\bf 52} (1985), 973-1007.

\bibitem{ShJ1974} J. A. Shalika, \emph{The multiplicitty one theorm for ${\rm GL}_n$}, Ann. Math. {\bf 100} (1974), 171-193.

\bibitem{Shin2011} S. W. Shin, \emph{Galois representations arising from some compact Shimura varieties}, Ann. Math. {\bf 173} (2011), 1645-1741.

\bibitem{Si1978}
A. Silberger, {\it The Langlands quotient theorem for $p$-adic groups}, Math. Ann. {\bf 236} (1978), 95-104.

\bibitem{SpBook} T. A. Springer, \emph{Algebraic groups}.

\bibitem{St2008} S. Stevens, \emph{The supercuspidal representations of $p$-adic classical groups}, Invent. Math. {\bf 172} (2008), 289-352.

\bibitem{TaS2014} S. Takeda, \emph{The twisted symmetric square $L$-function of ${\rm GL}(r)$}, Duke Math. J. {\bf 163} (2014), 175-266.

\bibitem{TaS2015} \rule{3em}{.5pt}\thinspace, \emph{On a certain metaplectic Eisenstein series and the twisted symmetric square $L$-function}, Math. Zeit. {\bf 281} (2015), 103-157.

\bibitem{Ta1950} J. T. Tate, \emph{Fourier analysis in number fields and Hecke's zeta-functions}, thesis, Princeton Univ., (1950) (published in [J. W. S. Cassels and A. Fr\"ohlich, editors, Algebraic number theory, $2^{\rm nd}$ ed., London Mathematical Society, 2010, pp. 305-347]).

\bibitem{Ta1979} J. T. Tate, \emph{Number Theoretic Background}, Proc. Symp. Pure Math. {\bf 33} (1979), vol. 2, 3-26.

\bibitem{Ti1979} J. Tits, \emph{Reductive groups over local fields}, Proc. Symp. Pure Math. {\bf 33} (1979), vol. 1, 29-69.

\bibitem{Tu1978} J. B. Tunnell, \emph{On the local Langlands conjecture for ${\rm GL}(2)$}, Invent. Math. {\bf 46} (1978), 179-200.

\bibitem{Wa2003} J.-L. Waldspurger, \emph{La formule de Plancherel pour les groupes $\mathfrak{p}$-adiques: d'apr\`es Harish-Chandra}, J. Inst. Math. Jussieu {\bf  2} (2003), 235-333.

\bibitem{Yu2001} J.-K. Yu, \emph{Construction of tame supercuspidal representations}, J.A.M.S. {\bf 14} (2001), 579-622.

\bibitem{We1951} A. Weil, \emph{Sur la Th\'eorie du Corps de Classes}, J. Math. Soc. Japan {\bf 3} (1951), 1-35.

\bibitem{Ze1980} A. V. Zelevinsky, \emph{Induced representations of reductive $\mathfrak{p}$-adic groups II. On irreducible representations of ${\rm GL}(n)$}, Ann. Sci. \'Ecole Norm. Sup. $4^{e}$ s\'erie {\bf 13} (1980), 165-210.
\end{thebibliography}
\end{document}